\documentclass[12pt]{article}

\usepackage{amssymb, amsmath, amsthm, amscd, graphicx}
\newtheorem{sch}{Scholium}[section]
\newtheorem{lem}{Lemma}[section]
\newtheorem{pro}{Proposition}[section]
\newtheorem{thm}{Theorem}[section]
\newtheorem{cor}{Corollary}[section]
\newtheorem{defn}{Definition}[section]
\newtheorem{rem}{Remark}[section]

\def\be{\begin{equation}}\def\ee{\end{equation}}
\def\colon{\,{:}\;}
\def\prf{\medbreak\noindent{\bf Proof}:\enspace}
\def\qed{\hspace*{\fill}\hbox{\vrule height 7pt \kern-.3pt
     \vbox{\hrule width 7pt
     \kern6.6pt\hrule width 7pt }\kern-.3pt\vrule height 7pt
     }\par}

\def\cirk{{\scriptstyle\circ}}
\def\ra{\rightarrow}
\def\BSigma{\mathbf\Sigma}
\def\ab{{\alpha,\beta}}
\def\abp{{\alpha^\prime,\beta^\prime}}
\def\abb{{\bar{\alpha},\bar{\beta}}}
\def\p{\prime}
\def\ud#1{\text{\b{#1}}}
\def\oph{\overline{\varphi}}
\def\ophi{\overline{\varphi}_\infty}
\def\tA{{\tt A}}

\def\i{{\tt i}}
\def\tp{{\tt p}}
\def\s{{\tt s}}

\def\N{{\bf N}}\def\Z{{\bf Z}}\def\R{{\bf R}}\def\Q{{\bf Q}}
\def\cC{{\mathcal C}}

\def\cF{{\mathcal F}}\def\cG{{\mathcal G}}

\def\cH{{\mathcal H}}
\def\cL{{\mathcal L}}

\begin{document}

\title{Computation of Topological Entropy via $\varphi$-expansion,\\
an Inverse Problem for the Dynamical Systems $\beta
x+\alpha\mod1$}

\author{B. Faller\footnote{e-mail: bastien.faller@epfl.ch} and
C.-E. Pfister\footnote{e-mail: charles.pfister@epfl.ch}\\
EPF-L, Institut d'analyse et calcul scientifique\\
CH-1015 Lausanne, Switzerland}

\date{16.05.2008
\\
}
\maketitle

\begin{abstract}
We give an algorithm, based on the $\varphi$-expansion of Parry,
in order to compute the topological entropy of a class of shift
spaces.  The idea is to solve an inverse problem for the dynamical
systems $\beta x +\alpha\mod 1$. The first part is an exposition
of the $\varphi$-expansion applied to piecewise monotone dynamical
systems. We formulate  for the validity of the
$\varphi$-expansion, necessary and sufficient conditions, which
are different from those in Parry's paper \cite{P2}.
\end{abstract}


\newpage

\section{Introduction}\label{section1}
\setcounter{equation}{0}

In 1957 R\'enyi published his paper \cite{R} about
representations for real numbers by $f$-expansions, called
hereafter $\varphi$-expansions, which had tremendous impact in
Dynamical Systems Theory. The ideas of R\'enyi were further
developed by Parry in \cite{P1} and \cite{P2}. See also the book
of Schweiger \cite{Sch}. The first part of the paper, section
\ref{section2}, is an exposition of the theory of
$\varphi$-expansions  in the setting of piecewise monotone
dynamical systems. Although many of the results of section
\ref{section2} are known, for example see  \cite{Bo} chapter 9
for Theorem \ref{thm2.5}, we state necessary and sufficient
conditions for the validity of the $\varphi$-expansion, which are
different from those in Parry's paper \cite{P2}, Theorem
\ref{thm2.1bis} and Theorem \ref{thm2.1ter}.

We then use  $\varphi$-expansions to study two interesting and
related problems in sections \ref{section3} and \ref{section4}.
When one applies the method of section \ref{section2} to the
dynamical system $\beta x+\alpha\mod1$, one obtains a symbolic
shift which is entirely described by two strings $\ud{u}^\ab$ and
$\ud{v}^\ab$ of symbols in a finite alphabet
$\tA=\{0,\ldots,k-1\}$. The shift space is given by
\begin{equation}\label{1.1}
\BSigma(\ud{u}^\ab,\ud{v}^\ab)=\big\{\ud{x}\in\tA^{\Z_+}\colon
\ud{u}^\ab\preceq\sigma^n\ud{x}\preceq\ud{v}^\ab\;\,\forall n\geq
0 \big\}\,,
\end{equation}
where $\preceq$ is the lexicographic order and $\sigma$ the shift
map. The particular case $\alpha=0$ has been much studied from
many different viewpoints ($\beta$-shifts). For $\alpha\not=0$ the
structure of the shift space is richer. A natural problem is to
study all shift spaces $\Sigma(\ud{u},\ud{v})$ of the form
\eqref{1.1} when we replace $\ud{u}^\ab$ and $\ud{v}^\ab$ by a
pair of strings $\ud{u}$ and $\ud{v}$. In section \ref{section3}
we give an algorithm, Theorem \ref{thm3.1}, based on the
$\varphi$-expansion, which allows to compute the topological
entropy of shift spaces $\Sigma(\ud{u},\ud{v})$. One of the
essential tool is the follower-set graph associated to the shift
space. This graph is presented in details in subsection
\ref{subsectionfollower}. The algorithm is given in subsection
\ref{subsectionalgo} and the computations of the topological
entropy in subsection \ref{topological}. The basic idea of the
algorithm is to compute two real numbers $\bar{\alpha}$ and
$\bar{\beta}$, given the strings $\ud{u}$ and $\ud{v}$, and to
show that the shift space $\BSigma(\ud{u},\ud{v})$ is a
modification of the shift space $\Sigma(\ud{u}^\abb,\ud{v}^\abb)$
obtained from the dynamical system
$\bar{\beta}x+\bar{\alpha}\mod1$, and that the topological
entropies of the two shift spaces are the same. In the last
section we consider the following inverse problem for the
dynamical systems $\beta x+\alpha \mod1$: given $\ud{u}$
 and $\ud{v}$, find $\alpha$
and $\beta$ so that
$$
\ud{u}=\ud{u}^\ab\quad\text{and}\quad\ud{v}=\ud{v}^\ab\,.
$$
The solution of this problem is given in Theorems \ref{thm4.1}
and \ref{thm4.2} for all $\beta>1$.

\section{$\varphi$-expansion for piecewise monotone dynamical \\
systems}\label{section2}
\setcounter{equation}{0}

\subsection{Piecewise monotone dynamical
systems}\label{subsection2.1}

Let $X:=[0,1]$ (with the euclidean distance). We consider the case
of piecewise monotone dynamical systems of the following type. Let
$0=a_0<a_1<\cdots<a_k=1$ and $I_j:=(a_j,a_{j+1})$, $j\in\tA$. We
set  $\tA:=\{0,\ldots,k-1\}$, $k\geq 2$, and
$$
S_0:=X\backslash \bigcup_{j\in\tA}I_j\,.
$$
For each $j\in\tA$ let
$$
f_j:I_j\mapsto J_j:=f_j(I_j)\subset [0,1]
$$
be a strictly monotone continuous map. When necessary we also
denote by $f_j$ the continuous extension of the map on the closure
$\overline{I}_j$ of $I_j$. We define a map $T$ on $X\backslash
S_0$ by setting
$$
T(x):=f_j(x)\quad \text{if $x\in I_j$}\,.
$$
The map $T$ is left undefined on $S_0$. We also assume that

\begin{equation}\label{2.1}
\big(\bigcup_{i\in\tA}J_i\big)\cap I_j=I_j\quad\forall j\,.
\end{equation}

\noindent
We introduce sets  $X_j$, $S_j$, and $S$  by setting for $j\geq 1$
$$
X_0:=[0,1]\,,\quad X_j:=X_{j-1}\backslash S_{j-1}\,,\quad
S_j:=\{x\in X_j\colon T(x)\in S_{j-1}\}\,,\quad S:=\bigcup_{j\geq
0}S_j\,.
$$

\begin{lem}\label{lem2.1}
Under the condition \eqref{2.1}, $T^n(X_{n+1})= X_1$ and
$T(X\backslash S)=X\backslash S$. $X\backslash S$ is dense in $X$.
\end{lem}

\prf Condition \eqref{2.1} is equivalent to $T(X_1)\supset X_1$.
Since $X_2=X_1\backslash S_1$ and $S_1=\{x\in X_1\colon
T(x)\not\in X_1\}$, we have $T(X_2)=X_1$. Suppose that
$T^n(X_{n+1})=X_1$; we prove that $T^{n+1}(X_{n+2})=X_{1}$. One
has $X_{n+1}=X_{n+2}\cup S_{n+1}$ and
$$
X_1=T^n(X_{n+1})=T^n(X_{n+2})\cup T^n(S_{n+1})\,.
$$
Applying once more $T$,
$$
X_1\subset T(X_1)=T^{n+1}(X_{n+2})\cup T^{n+1}(S_{n+1})\,.
$$
$T^{n+1}$ is defined on $X_{n+1}$ and $S_{n+1}\subset X_{n+1}$.
$$
T^{n+1}S_{n+1}=\{x\in X_{n+1}\colon T^{n+1}(x)\in S_0\}= \{x\in
X_{n+1}\colon T^{n+1}(x)\not\in X_1\}\,.
$$
Hence $T^{n+1}(X_{n+2})=X_{1}$. Clearly $T(X\backslash S)\subset
X\backslash S$ and $T(S\backslash S_0)\subset S$. Since $X_1$ is
the disjoint union of  $X\backslash S$ and $S\backslash S_0$, and
$TX_1\supset X_1$, we have $T(X\backslash S)=X\backslash S$. The
sets $X\backslash S_k$ are open and dense in $X$. By Baire's
Theorem $X\backslash S=\bigcap_{k}(X\backslash S_k)$ is dense.
\qed

Let $\Z_+:=\{0,1,2,\ldots\}$ and $\tA^{\Z_+}$ be equipped with
the product topology.  Elements of $\tA^{\Z_+}$ are called {\sf
strings} and denoted by $\ud{x}=(x_0,x_1,\ldots)$. A finite string
$\ud{w}=(w_0,\cdots,w_{n-1})$, $w_j\in\tA$, is a {\sf word}; we
also use the notation $\ud{w}=w_0\cdots w_{n-1}$. The {\sf length
of $\ud{w}$} is $|\ud{w}|=n$. A {\sf $n$-word} is a word of
length $n$. There is a single word of length $0$, the {\sf
empty-word} $\epsilon$. The set of all words is $\tA^*$. The
shift-map $\sigma\colon \tA^{\Z_+}\ra\tA^{\Z_+}$ is defined by
$$
\sigma(\ud{x}):=(x_1,x_2,\ldots)\,.
$$
We define two operations $\tp$ and $\s$ on
$\tA^*\backslash\{\epsilon\}$,
\begin{align}
\tp{\ud{w}}&:=\begin{cases} w_0\cdots w_{n-2}& \text{if
$\ud{w}=w_0\cdots w_{n-1}$ and $n\geq 2$}\label{opp}\\
\epsilon & \text{if $\ud{w}=w_0$}
\end{cases}\\
\s{\ud{w}}&:=\begin{cases} w_1\cdots w_{n-1}& \text{if
$\ud{w}=w_0\cdots w_{n-1}$ and $n\geq 2$}\label{ops}\\
\epsilon & \text{if $\ud{w}=w_0$.}
\end{cases}
\end{align}
On $\tA^{\Z_+}$ we define a total order, denoted by $\prec$. We
set
$$
\delta(j):=\begin{cases} +1 &\text{if $f_j$ is increasing}\\
-1 & \text{if $f_j$ is decreasing,}
\end{cases}
$$
and for a word $\ud{w}$,
$$
\delta(\ud{w}):=\begin{cases} \delta(w_0)\cdots \delta(w_{n-1})
&\text{if
$\ud{w}=w_0\cdots w_{n-1}$}\\
1& \text{if $\ud{w}=\epsilon$.}
\end{cases}
$$
Let $\ud{x}^\prime\not =\ud{x}^{\prime\prime}$ belong to
$\tA^{\Z_+}$; define  $j$ as the smallest integer with
$x^\prime_j\not=x^{\prime\prime}_j$. By definition
$$
\ud{x}^\prime\prec\ud{x}^{\prime\prime}\iff
\begin{cases}
x^\prime_j<x^{\prime\prime}_j &\text{if  $\delta(x^\prime_0\cdots
x^\prime_{j-1})=1$}\\
x^\prime_j>x^{\prime\prime}_j &\text{if  $\delta(x^\prime_0\cdots
x^\prime_{j-1})=-1$}\,.
\end{cases}
$$
As usual $\ud{x}^\prime\preceq \ud{x}^{\prime\prime}$ if and only
if $\ud{x}^\prime\prec\ud{x}^{\prime\prime}$ or
$\ud{x}^\prime=\ud{x}^{\prime\prime}$. When all maps $f_j$ are
increasing this order is the lexicographic order.

\subsection{$\varphi$-expansion}\label{subsection2.2}

We give an alternative description of a  piecewise monotone
dynamical system  as in Parry's paper \cite{P2}.  In this
description, when all maps $f_j$ are increasing, one could use
instead of the intervals $I_j$  the intervals
$I^\prime_j:=[a_j,a_{j+1})$, $j\in\tA$. In that case $S_0=\{a_k\}$
and $S_j=\emptyset$ for all $j\geq 1$. This would correspond to
the setting of Parry's paper \cite{P2}.

We define a map $\varphi$ on the disjoint union
$$
{\rm dom}\varphi:=\bigcup_{j=0}^{k-1}j+J_j\subset \R\,,
$$
by setting
\begin{equation}\label{2.3}
\varphi(x):= f^{-1}_j(t) \quad \text{if  $x=j+t$ and $t\in
J_j$}\,.
\end{equation}
The map $\varphi$ is  continuous, injective with range $X_1$. On
$X_1$ the inverse map is
$$
\varphi^{-1}(x)=j+Tx\quad\text{if $x\in I_j$}\,.
$$
For each $j$, such that $f_j$ is increasing, we define
$\overline{\varphi}^j$ on $j+[0,1]$ (using the extension of $f_j$
to $[a_j,a_{j+1}]$) by
\begin{equation}\label{2.4}
\overline{\varphi}^j(x):=\begin{cases} a_j & \text{if $x=j+t$ and
$t\leq
f_j(a_j)$}\\
f^{-1}_j(t) & \text{if  $x=j+t$ and $t\in J_j$}\\
a_{j+1} & \text{if  $x=j+t$ and $f_j(a_{j+1})\leq t$.}
\end{cases}
\end{equation}
For each $j$, such that $f_j$ is decreasing, we define
$\overline{\varphi}^j$ on $j+[0,1]$ by
\begin{equation}\label{2.5}
\overline{\varphi}^j(x):=\begin{cases} a_{j+1} & \text{if $x=j+t$
and $t\leq
f_j(a_{j+1})$}\\
f^{-1}_j(t) & \text{if  $x=j+t$ and $t\in J_j$}\\
a_{j} & \text{if  $x=j+t$ and $f_j(a_{j})\leq t$.}
\end{cases}
\end{equation}

\noindent
It is convenient below to consider the family of maps
$\overline{\varphi}^j$ as a single map defined on $[0,k]$, which
is denoted by $\overline{\varphi}$. In order to avoid ambiguities
at integers, where the map may be multi-valued, we always write a
point of $[j,j+1]$ as $x=j+t$, $t\in[0,1]$, so that
$$
\overline{\varphi}(j+t)\equiv\overline{\varphi}(x):=\overline{\varphi}^j(t)\,.
$$
We define the {\sf coding map} $\i:X\backslash S\ra \tA^{\Z_+}$ by
$$
\i(x):=(\i_0(x),\i_1(x),\ldots)\quad\text{with $\i_n(x):=j$ if
$T^nx\in I_j$}\,.
$$
The {\sf $\varphi$-code} of $x\in X\backslash S$ is the string
$\i(x)$, and we set
$$
\Sigma=\{\ud{x}\in\tA^{\Z_+}\colon \text{$\ud{x}=\i(x)$ for some
$x\in X\backslash S$}\}\,.
$$
For $x\in X\backslash S$ and any $n\geq 0$,
\begin{equation}\label{2.7}
\varphi^{-1}(T^nx)=\i_n(x)+T^{n+1}x\quad\text{and}\quad
\i(T^nx)=\sigma^n\i(x) \,.
\end{equation}
Let $z_j\in\tA$, $1\leq j \leq n$, and $t\in [0,1]$; we set
$$
\overline{\varphi}_1(z_1+t):=\overline{\varphi}(z_1+t)
$$
and
\begin{equation}\label{formule}
\overline{\varphi}_n(z_1,\ldots,z_n+t):=
\overline{\varphi}_{n-1}(z_1,\ldots,z_{n-1}+\overline{\varphi}(z_n+t))\,.
\end{equation}

\noindent
For $n\geq 1$ and $m\geq 1$ we have
\begin{equation}\label{formulegenerale}
\overline{\varphi}_{n+m}(z_1,\ldots,z_{n+m}+t) =
\overline{\varphi}_n(z_1,\ldots,z_{n}+
\overline{\varphi}_m(z_{n+1},\ldots,z_{n+m}+t))\,.
\end{equation}

\noindent
The map $t\mapsto \overline{\varphi}_n(x_0,\ldots,x_{n-1}+t)$ is
increasing if $\delta(x_0\cdots x_{n-1})=1$ and decreasing if
$\delta(x_0\cdots x_{n-1})=-1$. We also write
$\overline{\varphi}_n(\ud{x})$ for
$\overline{\varphi}_n(x_0,\ldots,x_{n-1})$.

\begin{defn}\label{defn2.2}
The real number $s$ has a {\sf $\varphi$-expansion}
$\ud{x}\in\tA^{\Z_+}$ if the following limit exists,
$$
s=\lim_{n\ra\infty}\overline{\varphi}_n(\ud{x})\equiv
\overline{\varphi}\big(x_0+\overline{\varphi}(x_1+\ldots)\big)\equiv
\overline{\varphi}_\infty(\ud{x})\,.
$$
The {\sf $\varphi$-expansion is well-defined} if for all
$\ud{x}\in\tA^{\Z_+}$,
$\lim_{n\ra\infty}\overline{\varphi}_n(\ud{x})=
\overline{\varphi}_\infty(\ud{x})$ exists.\\
The {\sf $\varphi$-expansion is valid} if for all $x\in
X\backslash S$ the $\varphi$-code $\i(x)$ of $x$ is a {\sf
$\varphi$-expansion} of $x$.
\end{defn}

If the $\varphi$-expansion is valid, then for $x\in X\backslash
S$, using \eqref{formulegenerale}, \eqref{2.7} and the continuity
of the maps $\overline{\varphi}^j$,
\begin{align}\label{2.8}
x&=\lim_{n\ra\infty}
\overline{\varphi}_n\big(\i_0(x),\ldots,\i_{n-1}(x)\big)\nonumber\\
&= \lim_{m\ra\infty}
\overline{\varphi}_{n}\big(\i_0(x),\ldots,\i_{n-1}(x)+
\overline{\varphi}_{m}\big(\i_n(x),\ldots,\i_{n+m-1}(x)\big)\big)
\\
&= \overline{\varphi}_{n}\big(\i_0(x),\ldots,\i_{n-1}(x)+
\overline{\varphi}_\infty(\i(T^nx)\big)\,.\nonumber
\end{align}
The basic and elementary fact of the $\varphi$-expansion is
\begin{equation}\label{2.9}
\text{$a,b\in[0,1]$ and $x_0<x_0^\prime$} \implies
\overline{\varphi}(x_0+a)\leq \overline{\varphi}(x_0^\prime+b)\,.
\end{equation}

We begin with two lemmas on the $\varphi$-code  (for Lemma
\ref{lem2.4} see e.g. \cite{CoE}).

\begin{lem}\label{lem2.4}
The $\varphi$-code $\i$ is $\preceq$-order-preserving on
$X\backslash S$: $x\leq y$ implies $\i(x)\preceq\i(y)$.
\end{lem}

\prf Let $x<y$. Either $\i_0(x)<\i_0(y)$, or $\i_0(x)=\i_0(y)$;
in the latter case, the strict monotonicity of $f_{\i_0(x)}$
implies
\begin{align*}
\varphi^{-1}(x)&=\i_0(x) +T(x)< \varphi^{-1}(y)=\i_0(x)
+T(y)\quad\text{if
$\delta(\i_0(x))=+1$}\\
\varphi^{-1}(x)&=\i_0(x) +T(x)> \varphi^{-1}(y)=\i_0(x)
+T(y)\quad\text{if $\delta(\i_0(x))=-1$.}
\end{align*}
Repeating this argument we get  $\i(x)\preceq\i(y)$.  \qed

\begin{lem}\label{lem2.4bis}
The $\varphi$-code $\i$ is continuous\footnote{\label{f2.3} If we
use the intervals $I^\prime_j=[a_j,a_{j+1})$, then we have only
right-continuity} on $X\backslash S$.
\end{lem}

\prf Let $x\in X\backslash S$ and $\{x^n\}\subset X\backslash S$,
$\lim_n x^n=x$.  Let $x\in I_{j_0}$. For $n$ large enough $x^n\in
I_{j_0}$ and $\i_0(x^n)=\i_0(x)=j_0$. Let $j_1:=\i_1(x)$; we can
choose $n_1$ so large that for $n\geq n_1$ $Tx_{n}\in I_{j_1}$.
Hence $\i_0(x^n)=j_0$ and $\i_1(x^n)=j_1$ for all $n\geq n_1$. By
induction we can find an increasing sequence $\{n_m\}$ such that
$n\geq n_m$ implies $\i_j(x)=\i_j(x^{n})$ for all
$j=0,\ldots,m$.  \qed

The next lemmas give the essential properties of the map
$\overline{\varphi}_\infty$.

\begin{lem}\label{lem2.2}
Let $\ud{x}\in\tA^{\Z_+}$. Then there exist $y_\uparrow(\ud{x})$
and $y_\downarrow(\ud{x})$ in $[0,1]$, such that
$y_\uparrow(\ud{x})\leq y_\downarrow(\ud{x})$;
$y_\uparrow(\ud{x})$ and $y_\downarrow(\ud{x})$ are the only
possible cluster points of the sequence
$\{\overline{\varphi}_n(\ud{x})\}_n$.\\
Let $x\in X\backslash S$ and set $\ud{x}:=\i(x)$. Then
$$
a_j\leq y_\uparrow(\ud{x})\leq x\leq y_\downarrow(\ud{x})\leq
a_{j+1}\quad \text{if $x_0=j$}\,.
$$
If the $\varphi$-expansion is valid, then each $y\in X\backslash
S$ has a unique $\varphi$-expansion\footnote{\label{f2.4} If we
use the intervals $I^\prime_j=[a_j,a_{j+1})$, this statement is
not correct.},
$$
y=\overline{\varphi}_\infty(\ud{x})\in X\backslash S\iff
\ud{x}=\i(y)\,.
$$
\end{lem}

\prf Consider the map
$$
t\mapsto \overline{\varphi}_n(x_0,\ldots,x_{n-1}+t)\,.
$$
Suppose that $\delta(x_0\cdots x_{n-1})=-1$. Then it is
decreasing, and for any $m$
\begin{align*}
\overline{\varphi}_{n+m}(x_0,\ldots,x_{n+m-1}) &=
\overline{\varphi}_n(x_0,\ldots,x_{n-1}+
\overline{\varphi}_m(x_n,\ldots,x_{n+m-1}))\\
&\leq \overline{\varphi}_n(x_0,\ldots,x_{n-1})\,.
\end{align*}
In particular the subsequence $\{\overline{\varphi}_n(\ud{x})\}_n$
of all $n$ such that $\delta(x_0\cdots x_{n-1})=-1$ is decreasing
with limit\footnote{\label{f2.limit} If the subsequence is finite,
then $y_\downarrow(\ud{x})$ is the last point of the subsequence.}
$y_\downarrow(\ud{x})$. When there is no $n$ such that
$\delta(x_0\cdots x_{n-1})=-1$, we set
$y_\downarrow(\ud{x}):=a_{x_0+1}$. Similarly, the subsequence
$\{\overline{\varphi}_n(\ud{x})\}_n$ of all $n$ such that
$\delta(x_0\cdots x_{n-1})=1$ is increasing with limit
$y_\uparrow(\ud{x})\leq y_\downarrow(\ud{x})$. When there is no
$n$ such that $\delta(x_0\cdots x_{n-1})=1$, we set
$y_\uparrow(\ud{x}):=a_{x_0}$. Since any
$\overline{\varphi}_n(\ud{x})$ appears in one of these sequences,
there are at most two cluster points for
$\{\overline{\varphi}_n(\ud{x})\}_n$.

Let $x\in X\backslash S$; $x=\varphi(\varphi^{-1}(x))$ and by
\eqref{2.7}
\begin{align}\label{identity}
x&=\varphi(\i_0(x)+Tx)=
\varphi(\i_0(x)+\varphi(\varphi^{-1}(Tx)))=
\varphi(\i_0(x)+\varphi(\i_1(x)+T^2x))=\cdots \nonumber\\
&=\varphi\big(\i_0(x)+\varphi(\i_1(x)+\ldots
+\varphi(\i_{n-1}(x)+T^{n}x))\big)\,.
\end{align}
By monotonicity
\begin{equation}\label{id1}
\big(\text{$x\in X\backslash S$ and
$\delta(\i_0(x)\cdots\i_{n-1}(x))=-1$}\big)\implies
\overline{\varphi}_n(\i_0(x),\cdots,\i_{n-1}(x))\geq x\,,
\end{equation}
and
\begin{equation}\label{id2}
\big(\text{$x\in X\backslash S$ and
$\delta(\i_0(x)\cdots\i_{n-1}(x))=1$}\big)\implies
\overline{\varphi}_n(\i_0(x),\cdots,\i_{n-1}(x))\leq x\,.
\end{equation}
The inequalities of Lemma \ref{lem2.2} follow from \eqref{id1},
\eqref{id2} and $\overline{\varphi}(\i_0(x)+t)\in [a_{x_0},
a_{x_0+1}]$.

Suppose that the $\varphi$-expansion is valid and that
$\overline{\varphi}_\infty(\ud{x})=y\in X\backslash S$. We prove
that $\ud{x}=\i(y)$. By hypothesis $y\in I_{x_0}$; using
\eqref{2.8} and the fact that $I_{x_0}$ is open, we can write
$$
y=\overline{\varphi}\big(x_0+\overline{\varphi}(x_1
+\overline{\varphi}(x_{2}+\ldots))\big)
=\varphi\big(x_0+\overline{\varphi}(x_1
+\overline{\varphi}(x_{2}+\ldots))\big)\,.
$$
This implies that
$$
\varphi^{-1}(y)=\i_0(y)+Ty= x_0+ \overline{\varphi}(x_1
+\overline{\varphi}(x_{2}+\ldots))\,.
$$
Since $Ty \in X\backslash S$, we can iterate this argument. \qed

\begin{lem}\label{lem2.3}
Let $\ud{x},\ud{x}^\prime\in\tA^{\Z_+}$ and
$\ud{x}\preceq\ud{x}^\prime$. Then any cluster point of
$\{\overline{\varphi}_n(\ud{x})\}_n$ is smaller then any cluster
point of $\{\overline{\varphi}_n(\ud{x}^\prime)\}_n$. In
particular, if $\overline{\varphi}_\infty$ is well-defined on
$\tA^{\Z_+}$, then $\overline{\varphi}_\infty$ is
order-preserving.
\end{lem}

\prf Let $\ud{x}\prec\ud{x}^\prime$ with $x_k=x^\prime_k$,
$k=0,\ldots,m-1$ and $x_m\not=x^\prime_m$. We have
$$
\overline{\varphi}_{m+n}(\ud{x})=
\overline{\varphi}_m(x_0,\ldots,x_{m-1}+
\overline{\varphi}_n(\sigma^m\ud{x}))\,.
$$
By \eqref{2.9}, if $\delta(x_0\cdots x_{m-1})=1$, then
$x_m<x^\prime_m$ and for any $n\geq 1$, $\ell\geq 1$,
$$
\overline{\varphi}_n(\sigma^m\ud{x})=\overline{\varphi}_1(x_m+
\overline{\varphi}_{n-1}(\sigma^{m+1}\ud{x}))\leq
\overline{\varphi}_\ell(\sigma^m\ud{x}^\prime)=
\overline{\varphi}_1(x^\prime_m+
\overline{\varphi}_{\ell-1}(\sigma^{m+1}\ud{x}^\prime))\,;
$$
if $\delta(x_0\cdots x_{m-1})=-1$, then $x_m>x^\prime_m$ and
$$
\overline{\varphi}_n(\sigma^m\ud{x})=\overline{\varphi}_1(x_m+
\overline{\varphi}_{n-1}(\sigma^{m+1}\ud{x}))\geq
\overline{\varphi}_\ell(\sigma^m\ud{x}^\prime)=
\overline{\varphi}_1(x^\prime_m+
\overline{\varphi}_{\ell-1}(\sigma^{m+1}\ud{x}^\prime))\,.
$$
Therefore, in both cases, for any $n\geq 1$, $\ell\geq 1$,
$$
\overline{\varphi}_{m+n}(\ud{x})\leq
\overline{\varphi}_{m+\ell}(\ud{x}^\prime)\,.
$$
\qed

\begin{lem}\label{lem2.4ter}
Let $\ud{x}\in\tA^{\Z_+}$ and $x_0=j$.\\
1) Let $\delta(j)=1$ and $y_\uparrow(\ud{x})\in \overline{I}_j$ be
a cluster point of $\{\overline{\varphi}_n(\ud{x})\}$. Then
$f_j\big(y_\uparrow(\ud{x})\big)\geq y_\uparrow(\sigma\ud{x})$ if
$y_\uparrow(\ud{x})=a_j$, $f_j\big(y_\uparrow(\ud{x})\big)\leq
y_\uparrow(\sigma\ud{x})$ if $y_\uparrow(\ud{x})=a_{j+1}$ and
$f_j\big(y_\uparrow(\ud{x})\big)=y_\uparrow(\sigma\ud{x})$
otherwise. The same conclusions hold when $y_\downarrow(\ud{x})$
is a
cluster point of $\{\overline{\varphi}_n(\ud{x})\}$.\\
2) Let $\delta(j)=-1$ and $y_\uparrow(\ud{x})\in \overline{I}_j$
be a cluster point of $\{\overline{\varphi}_n(\ud{x})\}$. Then
$f_j\big(y_\uparrow(\ud{x})\big)\leq y_\downarrow(\sigma\ud{x})$
if $y_\uparrow(\ud{x})=a_j$, $f_j\big(y_\uparrow(\ud{x})\big)\geq
y_\downarrow(\sigma\ud{x})$ if $y_\uparrow(\ud{x})=a_{j+1}$ and
$f_j\big(y_\uparrow(\ud{x})\big)=y_\downarrow(\sigma\ud{x})$
otherwise. The same conclusions hold when  $y_\downarrow(\ud{x})$
is a cluster point of $\{\overline{\varphi}_n(\ud{x})\}$.
\end{lem}

\prf Set $f_j(\overline{I}_j):=[\alpha_j,\beta_j]$. Suppose for
example that $\delta(j)=-1$ and that $n_k$ is the subsequence of
all $m$ such that $\delta(x_0,\ldots,x_m)=1$. Since $\delta(j)=-1$
the sequence $\{\oph_{n_k-1}(\sigma\ud{x})\}_k$ is decreasing.
Hence by continuity
\begin{equation}\label{mon}
y_\uparrow(\ud{x})=\lim_k\overline{\varphi}_{n_k}(\ud{x})=
\overline{\varphi}(j+\lim_k\overline{\varphi}_{n_k-1}(\sigma\ud{x}))=
\overline{\varphi}(j+y_\downarrow(\sigma\ud{x}))\,.
\end{equation}
If $y_\uparrow(\ud{x})=a_j$, then $f_j(a_j)=\beta_j\leq
y_\downarrow(\sigma\ud{x})$; if $y_\uparrow(\ud{x})=a_{j+1}$, then
$f_j(a_{j+1})=\alpha_j\geq y_\downarrow(\sigma\ud{x})$; if
$a_j<y_\uparrow(\ud{x})<a_{j+1}$, then
$$
j+f_j\big(y_\uparrow(\ud{x})\big)=
\varphi^{-1}\big(\varphi(j+\lim_k\overline{\varphi}_{n_k-1}(\sigma\ud{x}))\big)=
j+y_\downarrow(\sigma\ud{x})\,.
$$
Similar proofs for the other cases. \qed

\begin{lem}\label{lem2.4quatro}
Let $\ud{x}\in\tA^{\Z_+}$.\\
1) If $\{\overline{\varphi}_n(\ud{x})\}$ has two cluster points,
and if  $y\in \big(y_\uparrow(\ud{x}),
y_\downarrow(\ud{x})\big)$, then
$y\in X\backslash S$, $\i(y)=\ud{x}$ and $y$ has no $\varphi$-expansion.\\
Let $x\in X\backslash S$ and set $\ud{x}:=\i(x)$.\\
2) If $\lim_n\overline{\varphi}_n(\ud{x})=y_\uparrow(\ud{x})$ and
if $y\in (y_\uparrow(\ud{x}),
x)$, then $y\in X\backslash S$, $\i(y)=\ud{x}$ and  $y$ has no $\varphi$-expansion.\\
3) If $\lim_n\overline{\varphi}_n(\ud{x})=y_\downarrow(\ud{x})$
and if $y\in(x,y_\downarrow(\ud{x}))$, then $y\in X\backslash S$,
$\i(y)=\ud{x}$ and $y$ has no $\varphi$-expansion.
\end{lem}

\prf Suppose that $y_\uparrow(\ud{x})<y< y_\downarrow(\ud{x})$.
Then $y\in I_{x_0}$ and $\i_0(y)=x_0$. From Lemma \ref{lem2.4ter}
$$
y_\uparrow(\sigma\ud{x})<Ty<y_\downarrow(\sigma\ud{x})\quad
\text{if $\delta(x_0)=1$}\,,
$$
and
$$
y_\downarrow(\sigma\ud{x})>Ty>y_\uparrow(\sigma\ud{x})\quad
\text{if $\delta(x_0)=-1$}\,.
$$
Iterating this argument we prove that $T^ny\in I_{x_n}$ and
$\i_n(y)=x_n$ for all $n\geq 1$. Suppose that $y$ has a
$\varphi$-expansion, $y=\overline{\varphi}_\infty(\ud{x}^\prime)$.
If $\ud{x}^\p\prec\ud{x}$, then by Lemma \ref{lem2.3}
$\ophi(\ud{x}^\p)\leq y_\uparrow(\ud{x})$ and if $\ud{x}\prec
\ud{x}^\p$, then by Lemma \ref{lem2.3} $y_\downarrow(\ud{x})\leq
\ophi(\ud{x}^\p)$, which leads to a contradiction. Similar proofs
in cases 2 and 3. \qed

\begin{lem}\label{lem2.4quinto}
Let $\ud{x}^\prime\in\tA^{\Z_+}$ and $x\in X\backslash S$. Then
$$
\text{$y_\downarrow(\ud{x}^\prime)<x$ $\implies$
$\ud{x}^\prime\preceq \i(x)$}
\quad\text{and}\quad\text{$x<y_\uparrow(\ud{x}^\prime)$ $\implies$
$\i(x)\preceq \ud{x}^\prime$.}
$$
\end{lem}

\prf Suppose that $y_\downarrow(\ud{x}^\prime)<x$ and
$y_\downarrow(\ud{x}^\prime)$ is a cluster point. Either
$x^\prime_0<\i_0(x)$ or $x^\prime_0=\i_0(x)$ and by Lemma
\ref{lem2.4ter}
$$
y_\downarrow(\sigma\ud{x}^\prime) <Tx\quad\text{if
$\delta(x_0^\prime)=1$,}
$$
or
$$
y_\uparrow(\sigma\ud{x}^\prime)
>Tx\quad\text{if $\delta(x_0^\prime)=-1$.}
$$
Since  $y_\downarrow(\sigma\ud{x}^\prime)$ or
$y_\uparrow(\sigma\ud{x}^\prime)$ is a cluster point we can repeat
the argument and conclude that $\ud{x}^\prime\preceq\i(x)$. If
$y_\downarrow(\ud{x}^\prime)$ is not a cluster point, then we use
the cluster point
$y_\uparrow(\ud{x}^\prime)<y_\downarrow(\ud{x}^\prime)$ for the
argument. \qed

\begin{thm}\label{thm2.1}{\rm \cite{P2}}
A $\varphi$-expansion is valid if and only if the $\varphi$-code
$\i$ is injective on $X\backslash S$.
\end{thm}

\prf Suppose that the $\varphi$-expansion is valid. If $x\not=z$,
then
$$
x=\overline{\varphi}\big(\i_0(x)+
\overline{\varphi}(\i_1(x)+\ldots)\big)\not=
\overline{\varphi}\big(\i_0(z)+
\overline{\varphi}(\i_1(z)+\ldots)\big)=z\,,
$$
and therefore $\i(x)\not=\i(z)$. Conversely, assume  that
$x\not=z$ implies $\i(x)\not=\i(z)$. Let $x\in X\backslash S$,
$\ud{x}=\i(x)$, and suppose for example that $y_\uparrow(\ud{x})<
y_\downarrow(\ud{x})$ are two cluster points. Then by Lemma
\ref{lem2.4quatro} any $y$ such that $y_\uparrow(\ud{x})<y<
y_\downarrow(\ud{x})$ is in $X\backslash S$ and $\i(y)=\ud{x}$,
contradicting the hypothesis. Therefore
$z:=\lim_n\overline{\varphi}_n(\ud{x})$ exists. If $z\not=x$,
then we get again a contradiction using Lemma \ref{lem2.4quatro}.
\qed

Theorem \ref{thm2.1} states that the validity of the
$\varphi$-expansion is equivalent to the injectivity of the map
$\i$ defined on $X\backslash S$. One can also state that the
validity of the $\varphi$-expansion is equivalent to the
surjectivity  of the map $\overline{\varphi}_\infty$.

\begin{thm}\label{thm2.1bis}
A $\varphi$-expansion is valid if and only if
$\overline{\varphi}_\infty\colon \tA^{\Z_+}\ra [0,1]$ is
well-defined on $\tA^{\Z_+}$ and surjective.
\end{thm}

\prf Suppose that the $\varphi$-expansion is valid. Let
$\ud{x}\in\tA^{\Z_+}$ and suppose that
$\{\overline{\varphi}_n(\ud{x})\}_n$ has two different
accumulation points $y_\uparrow<y_\downarrow$. By Lemma
\ref{lem2.4quatro} we get a contradiction. Thus
$\overline{\varphi}_\infty(\ud{x})$ is well-defined for any
$\ud{x}\in\tA^{\Z_+}$.

To prove the surjectivity of $\overline{\varphi}_\infty$ it is
sufficient to consider $s\in S$. The argument is a variant of the
proof of Lemma \ref{lem2.4quatro}. Let $\ud{x}^\prime$ be a string
such that for any $n\geq 1$
$$
f_{x^\prime_{n-1}}\cirk \cdots \cirk
f_{x^\prime_0}(s)\in\overline{I}_{x^\prime_n}\,.
$$
We use here the extension of $f_j$ to $\overline{I}_j$; we have a
choice for $x^\prime_n$ whenever $f_{x^\prime_{n-1}}\cirk \cdots
\cirk f_{x^\prime_0}(s)\in S_0$. Suppose that
$\overline{\varphi}_\infty(\ud{x}^\prime)<s$ and that
$\overline{\varphi}_\infty(\ud{x}^\prime)<z<s$. Since $s,
\overline{\varphi}_\infty(\ud{x}^\prime)\in
\overline{I}_{x_0^\prime}$, we have $z\in I_{x_0^\prime}$ and
therefore $\i(z)=x_0^\prime$. Moreover,
$$
\overline{\varphi}_\infty(\sigma\ud{x}^\prime)<Tz<f_{x^\prime_0}(s)\quad\text{if
$\delta(x_0^\prime)=1$}
$$
or
$$
f_{x^\prime_0}(s)<Tz<\overline{\varphi}_\infty(\sigma\ud{x}^\prime)\quad\text{if
$\delta(x_0^\prime)=-1$}\,.
$$
Iterating the argument we get $z\in X\backslash S$ and
$\i(z)=\ud{x}^\prime$, contradicting the validity of the
$\varphi$-expansion. Similarly  we exclude the possibility that
$\overline{\varphi}_\infty(\ud{x}^\prime)>s$, thus proving the
surjectivity of the map $\overline{\varphi}_\infty$.

Suppose that $\overline{\varphi}_\infty\colon \tA^{\Z_+}\ra [0,1]$
is well-defined and surjective. Let $x\in X\backslash S$ and
$\ud{x}=\i(x)$. Suppose that $x<
\overline{\varphi}_\infty(\ud{x})$. By Lemma \ref{lem2.4quatro}
any $z$, such that $x<z<\overline{\varphi}_\infty(\ud{x})$,
 does not have a $\varphi$-expansion. This contradicts
the hypothesis that $\overline{\varphi}_\infty$ is surjective.
Similarly we exclude the possibility that $x>
\overline{\varphi}_\infty(\ud{x})$. \qed

\begin{thm}\label{thm2.1ter}
A $\varphi$-expansion is valid if and only if
$\overline{\varphi}_\infty: \tA^{\Z_+}\ra [0,1]$ is well-defined,
continuous and there exist $\ud{x}^+$ with
$\overline{\varphi}_\infty(\ud{x}^+)=1$ and $\ud{x}^-$ with
$\overline{\varphi}_\infty (\ud{x}^-)=0$.
\end{thm}

\prf Suppose that the $\varphi$-expansion is valid. By Theorem
\ref{thm2.1bis} $\overline{\varphi}_\infty$ is well-defined and
surjective so that there exist $\ud{x}^+$ and $\ud{x}^-$ with
$\overline{\varphi}_\infty(\ud{x}^+)=1$ and
$\overline{\varphi}_\infty(\ud{x}^-)=0$. Suppose that
$\ud{x}^n\downarrow\ud{x}$ and set
$y:=\overline{\varphi}_\infty(\ud{x})$,
$x_n:=\overline{\varphi}_\infty(\ud{x}^n)$. By Lemma \ref{lem2.3}
the sequence $\{x_n\}$ is monotone decreasing; let $x:=\lim_nx_n$.
Suppose that  $y<x$ and $y<z<x$. Since $y<z<x_n$ for any $n\geq
1$ and $\lim_n\ud{x}^n=\ud{x}$, we prove, as in the beginning of
the proof of Lemma \ref{lem2.4quatro}, that $z\in X\backslash S$.
The validity of the $\varphi$-expansion implies that
$z=\overline{\varphi}_\infty(\i(z))$. By Lemma \ref{lem2.4quinto}
$$
\ud{x}\preceq\i(z)\preceq\ud{x}^n\,.
$$
Since these inequalities are valid for any $z$, with $y<z<x$, the
validity of $\varphi$-expansion implies that we have strict
inequalities, $\ud{x}\prec\i(z)\prec \ud{x}^n$. This contradicts
the hypothesis that $\lim_{n\ra\infty}\ud{x}^n=\ud{x}$. A similar
argument holds in the case $\ud{x}^n\uparrow\ud{x}$. Hence
$$
\lim_{n\ra\infty}\ud{x}^n=\ud{x}\implies \lim_{n\ra\infty}
\overline{\varphi}_\infty(\ud{x}^n)=\overline{\varphi}_\infty(\ud{x})\,.
$$

Conversely, suppose that $\overline{\varphi}_\infty \colon
\tA^{\Z_+}\ra [0,1]$ is well-defined and continuous. Then, given
$\delta>0$ and $\ud{x}\in\tA^{\Z_+}$, $\exists n$ so that
$$
0\leq\sup\{\overline{\varphi}_\infty(\ud{x}^\prime)\colon
x^\prime_j=x_j\;j=0,\ldots,n-1\}-
\inf\{\overline{\varphi}_\infty(\ud{x}^\prime)\colon
x^\prime_j=x_j\;j=0,\ldots,n-1\}\leq \delta\,.
$$
We set
$$
\ud{x}^{n,-}:=x_0\cdots x_{n-1}\ud{x}^-\quad\text{and}\quad
\ud{x}^{n,+}:=x_0\cdots x_{n-1}\ud{x}^+\,.
$$
For any $x\in X\backslash S$ we have the identity
\eqref{identity},
$$
x=\varphi\big(\i_0(x)+\varphi(\i_1(x)+\ldots
+\varphi(\i_{n-1}+T^{n}x))\big)=
\overline{\varphi}_{n}(\i_0(x),\ldots,\i_{n-1}(x)+T^{n}x)\,.
$$
If $\delta(\i_0(x)\cdots\i_{n-1}(x))=1$, then
\begin{align*}
\overline{\varphi}_\infty(\ud{x}^{n,-}):&= \overline{\varphi}_n(
\i_0(x),\ldots,\i_{n-1}(x)+\overline{\varphi}_\infty(\ud{x}^-))\\
&=\overline{\varphi}_n(
\i_0(x),\ldots,\i_{n-1}(x))\\
&\leq
\overline{\varphi}_{n}(\i_0(x),\ldots,\i_{n-1}(x)+T^{n}x)\\
&\leq \overline{\varphi}_{n}(
\i_0(x),\ldots,\i_{n-1}(x)+1)\\
&=\overline{\varphi}_{n}(
\i_0(x),\ldots,\i_{n-1}(x)+\overline{\varphi}_\infty(\ud{x}^+))=:
\overline{\varphi}_\infty(\ud{x}^{n,+})\,.
\end{align*}
If $\delta(\i_0(x)\cdots\i_{n-1}(x))=-1$, then the inequalities
are reversed. Letting $n$ going to infinity, we get
$\overline{\varphi}_\infty(\i(x))=x$. \qed

\begin{rem}\label{rem2.4}
When the maps  $f_0$  and $f_{k-1}$ are  increasing, then we can
take
$$
\ud{x}^+=(k-1,k-1,\ldots)\quad\text{and}\quad
\ud{x}^{-}=(0,0,\ldots)\,.
$$
\end{rem}

\begin{thm}\label{thm2.2}{\rm \cite{P2}}
A necessary and sufficient condition for a $\varphi$-expansion to
be valid is that $S$ is dense in $[0,1]$. A sufficient condition
is $\sup_t|\varphi^\prime(t)|<1$.
\end{thm}

\bigskip

For each $j\in\tA$ we define (the limits are taken with $x\in
X\backslash S$)
\begin{equation}\label{2.12}
\ud{u}^j:=\lim_{x\downarrow a_{j}}\i(x)\quad\text{and}\quad
\ud{v}^j :=\lim_{x\uparrow a_{j+1}}\i(x)\,.
\end{equation}
The strings $\ud{u}^j$ and $\ud{v}^j$ are called {\sf virtual
itineraries}. Notice that $\underline{v}^j
\prec\underline{u}^{j+1}$ since $v^j_0<u_0^{j+1}$.
\begin{equation}\label{2.12b}
\sigma^k\ud{u}^j=\sigma^k(\lim_{x\downarrow a_{j}}\i(x))=
\lim_{x\downarrow a_{j}}\sigma^k\i(x)=\lim_{x\downarrow
a_{j}}\i(T^kx)\qquad(x\in X\backslash S)\,.
\end{equation}

\begin{pro}\label{prof2.2}
Suppose that $\ud{x}^\prime\in\tA^{\Z_+}$ verifies
$\ud{u}^{x_n^\prime}\prec \sigma^n\ud{x}^\prime\prec
\ud{v}^{x_n^\prime}$ for all $n\geq 0$. Then there exists $x\in
X\backslash S$ such that $\i(x)=\ud{x}^\prime$.
\end{pro}

\noindent
Notice that we do not assume that the $\varphi$-expansion is valid
or that the map $\overline{\varphi}_\infty$ is well-defined. For
unimodal maps see e.g. Theorem II.3.8 in \cite{CoE}. Our proof is
different.

\prf If $y_\uparrow(\ud{x}^\prime)<y_\downarrow(\ud{x}^\prime)$
are two cluster points, then this follows from Lemma
\ref{lem2.4quatro}. Therefore, assume that
$\lim_{n}\overline{\varphi}_n(\ud{x}^\prime)$ exists. Either
there exists $m>1$ so that
$y_\uparrow(\sigma^m\ud{x}^\prime)<y_\downarrow(\sigma^m\ud{x}^\prime)$
are two cluster points, or
$\lim_{n}\overline{\varphi}_n(\sigma^m\ud{x}^\prime)$ exists for
all $m\geq 1$.

In the first case, there exists $z_m\in X\backslash S$,
$$
y_\uparrow(\sigma^m\ud{x}^\prime)<z_m<y_\downarrow(\sigma^m\ud{x}^\prime)
\quad\text{and}\quad \i(z_m)=\sigma^m\ud{x}^\prime\,.
$$
Let
$$
z_{m-1}:=\overline{\varphi}(x^\prime_{m-1}+z_m)\,.
$$
We show that  $a_{x^\prime_{m-1}}< z_{m-1}<a_{x^\prime_{m-1}+1}$.
This implies that $z_m\in {\rm int(dom} \varphi)$ so that
$$
\varphi^{-1}(z_{m-1})=x^\prime_{m-1}+Tz_{m-1}=x^\prime_{m-1}+z_m\,.
$$
Suppose that $\delta(x^\prime_{m-1})=1$ and
$a_{x^\prime_{m-1}}=z_{m-1}$. Then for any $y\in X\backslash S$,
$y>a_{x^\prime_{m-1}}$, we have $Ty>z_{m}$. Therefore, by Lemma
\ref{lem2.4}, $\i(Ty)\succeq \i(z_m)=\sigma^m\ud{x}^\prime$;
$\i_0(y)=x^\prime_{m-1}$ when $y$ is close to
$a_{x^\prime_{m-1}}$, so that
$$
\lim_{y\downarrow
a_{x^\prime_{m-1}}}\i(y)=\ud{u}^{x^\prime_{m-1}}\succeq
\sigma^{m-1}\ud{x}^\prime\,,
$$
which is a contradiction.  Similarly we exclude the cases
$\delta(x^\prime_{m-1})=1$ and $a_{x^\prime_{m-1}+1}=z_{m-1}$,
$\delta(x^\prime_{m-1})=-1$ and $a_{x^\prime_{m-1}}=z_{m-1}$,
$\delta(x^\prime_{m-1})=-1$ and $a_{x^\prime_{m-1}+1}=z_{m-1}$.
Iterating this argument we get the existence of $z_0\in
X\backslash S$ with $\i(z_0)=\ud{x}^\prime$.

In the second case,
$\lim_{n}\overline{\varphi}_n(\sigma^m\ud{x}^\prime)$ exists for
all $m\geq 1$. Let
$x:=\lim_{n}\overline{\varphi}_n(\ud{x}^\prime)$. Suppose that
$x^\prime_0=j$, so that $\ud{u}^{j}\prec\ud{x}^\prime\prec
\ud{v}^{j}$. By Lemma \ref{lem2.4} and definition of $\ud{u}^{j}$
and $\ud{v}^{j}$ there exist $z_1,z_2\in I_j$ such that
$$
z_1<x<z_2\quad\text{and}\quad\ud{u}^{j}\preceq\i(z_1)\prec
\ud{x}^\prime\prec\i(z_2)\preceq \ud{v}^{j}\,.
$$
Therefore $a_j< x<a_{j+1}$, $\i_0(x)=x_0^\prime$ and $Tx=
\overline{\varphi}_\infty(\sigma\ud{x}^\prime)$ (Lemma
\ref{lem2.4ter}). Iterating  this argument we get
$\ud{x}^\prime=\i(x)$.  \qed

\begin{thm}\label{thm2.5}
Suppose that the $\varphi$-expansion is valid. Then
\begin{enumerate}
\item
$\Sigma:=\{\i(x)\in\tA^{\Z_+}\colon x\in X\backslash
S\}=\{\ud{x}\in\tA^{\Z_+}\colon \ud{u}^{x_n}\prec
\sigma^n\ud{x}\prec \ud{v}^{x_n}\quad\forall\,n\geq 0\}$.
\item
The map $\i\colon X\backslash S\ra \Sigma$ is bijective,
$\overline{\varphi}_\infty\cirk\i={\rm id}$ and
$\i\cirk\overline{\varphi}_\infty={\rm id}$.\\
Both maps $\i$ and $\overline{\varphi}_\infty$ are
order-preserving.
\item
$\sigma(\Sigma)=\Sigma$ and
$\overline{\varphi}_\infty(\sigma\ud{x})=T\overline{\varphi}_\infty(\ud{x})$
if $\ud{x}\in\Sigma$.
\item
If $\ud{x}\in\tA^{\Z_+}\backslash\Sigma$, then there exist
$m\in\Z_+$ and $j\in\tA$ such that
$\overline{\varphi}_\infty(\sigma^m\ud{x})=a_j$.
\item
$\forall n\geq 0\,,\,\forall j\in\tA\colon
\quad\ud{u}^{u^j_n}\preceq \sigma^n\ud{u}^j\prec \ud{v}^{u^j_n}$
if $\delta(\ud{u}^j_0\cdots\ud{u}^j_{n-1})=1$ and
$\ud{u}^{u^j_n}\prec \sigma^n\ud{u}^j\preceq \ud{v}^{v^j_n}$ if
$\delta(\ud{u}^j_0\cdots\ud{u}^j_{n-1})=-1$.
\item
$\forall n\geq 0\,,\,\forall j\in\tA\colon
\quad\ud{u}^{u^j_n}\preceq \sigma^n\ud{v}^j\prec \ud{v}^{u^j_n}$
if $\delta(\ud{v}^j_0\cdots\ud{v}^j_{n-1})=-1$ and
$\ud{u}^{u^j_n}\prec \sigma^n\ud{v}^j\preceq \ud{v}^{v^j_n}$ if
$\delta(\ud{u}^j_0\cdots\ud{u}^j_{n-1})=1$.
\end{enumerate}
\end{thm}

\prf Let $x\in X\backslash S$. Clearly, by monotonicity,
$$
\ud{u}^{\i_k(x)}\preceq \sigma^k\i(x)\preceq \ud{v}^{\i_k(x)}\quad
\forall\;k\in\Z_+\,.
$$
Suppose that  there exist $x\in X\backslash S$ and $k$ such that
 $\sigma^k\i(x)= \ud{v}^{\i_k(x)}$. Since
$(\sigma^k\i(x))_0=\i_0(T^kx)$, we can assume, without
restricting the generality, that $k=0$ and $\i_0(x)=j$. Therefore
$x\in (a_{j}, a_{j+1})$, and for all $y\in X\backslash S$, such
that $x\leq y<a_{j+1}$, we have by Lemma \ref{lem2.4} that
$\i(y)=\i(x)=\ud{v}^{j}$. By Theorem \ref{thm2.1} this
contradicts the hypothesis that the $\varphi$-expansion is valid.
The other case, $\sigma^k\i(x)= \ud{u}^{\i_k(x)}$, is treated
similarly. This proves half of the first statement. The second
half is a consequence of Proposition \ref{prof2.2}. The second
statement also follows, as well as the third, since
$T(X\backslash S)=X\backslash S$ (we assume that \eqref{2.1}
holds).

Let $\ud{x}\in\tA^{\Z_+}\backslash\Sigma$ and $m\in\Z_+$ be the
smallest integer such that one of the conditions defining $\Sigma$
is not verified. Then either
$\sigma^m\ud{x}\preceq\underline{u}^{x_m}$, or
$\sigma^m\ud{x}\succeq \underline{v}^{x_m}$. The map
$\overline{\varphi}_\infty$ is continuous (Theorem
\ref{thm2.1ter}). Hence, for any $j\in\tA$,
$$
\overline{\varphi}_\infty(\ud{u}^j)=a_{j}\quad\text{and}\quad
\overline{\varphi}_\infty(\ud{v}^j)=a_{j+1}\,.
$$
Let $\sigma^m\ud{x}\preceq\underline{u}^{x_m}$. Since
$\ud{v}^{x_m-1}\prec\sigma^m\ud{x}$,
$$
a_{x_m}=\overline{\varphi}_\infty(\ud{v}^{x_m-1})
\leq\overline{\varphi}_\infty(\sigma^m\ud{x})\leq
\overline{\varphi}_\infty(\ud{u}^{x_m})=a_{x_m}\,.
$$
The other case is treated in the same way. From definition
\eqref{2.12} $\ud{u}^{u^j_n}\preceq \sigma^n\ud{u}^j\preceq
\ud{v}^{u^j_n}$. Suppose that
$\delta(\ud{u}^j_0\cdots\ud{u}^j_{n-1})=1$ and
$\sigma^n\ud{u}^j=\ud{v}^{u^j_n}$. By continuity of the
$\varphi$-code there exists $x\in X\backslash S$ such that
$x>a_j$ and $\i_k(x)=\ud{u}^j_k$, $k=0,\ldots,n$. Let $a_j<y<x$.
Since $\delta(\ud{u}^j_0\cdots\ud{u}^j_{n-1})=1$, $T^ny<T^nx$ and
consequently
$$
\lim_{y\downarrow
a_j}\i(T^ny)=\sigma^n\ud{u}^j\preceq\i(T^nx)\preceq
\ud{v}^{u^j_n}\,.
$$
Hence $\sigma^n\i(x)= \ud{v}^{x_n}$, which is a contradiction. The
other cases are treated similarly. \qed

\subsection{Dynamical system $\beta x+\alpha\mod 1$}\label{subsection2.3}

We consider the family of dynamical systems $\beta x+\alpha\mod 1$
with $\beta >1$ and $0\leq \alpha<1$. For given $\alpha$ and
$\beta$, the dynamical system is described by $k=\lceil
\alpha+\beta\rceil$ intervals $I_j$ and maps $f_j$,
$$
I_0=\Big(0,\frac{1-\alpha}{\beta}\Big)\,,\,I_j=\Big(\frac{j-\alpha}{\beta},
\frac{j+1-\alpha}{\beta}\Big)
\,,\,j=1,\ldots,k-2\,,\,I_{k-1}=\Big(\frac{k-1-\alpha}{\beta},1\Big)
$$
and
$$
f_j(x)=\beta x+\alpha-j\,,\,j=0,\ldots k-1\,.
$$
The maps $T_{\alpha,\beta}$, $\varphi^\ab$ and
$\overline{\varphi}^\ab$ are defined as in subsection
\ref{subsection2.1}.
$$
\overline{\varphi}^\ab(x)=\begin{cases} 0 & \text{if
$0\leq x\leq\alpha$}\\
\beta^{-1}(x-\alpha) & \text{if $\alpha\leq x\leq \alpha+\beta$} \\
1 & \text{if $\alpha+\beta\leq x\leq \lceil \alpha+\beta\rceil$}
\end{cases}
$$
and
\begin{equation}\label{eqso}
S_0=\{a_j\colon j=1,\ldots,k-1\}\cup\{0,1\}\quad\text{with}\quad
a_j:=\beta^{-1}(j-\alpha)\,.
\end{equation}
Since all maps are increasing the total order on $\tA^{\Z_+}$ is
the lexicographic order. We have $2k$ virtual orbits, but only two
of them are important. Indeed, if we set
$$
\ud{u}^{\alpha,\beta}:=\ud{u}^0\quad\text{and}
\quad\ud{v}^{\alpha,\beta}:=\ud{v}^{k-1}\,,
$$
then
$$
\ud{u}^j=j\ud{u}^{\alpha,\beta}\,,\,j=1,\ldots k-1
$$
and
$$
\ud{v}^j=j\ud{v}^{\alpha,\beta}\,,\,j=0,\ldots,k-2\,.
$$

\begin{pro}\label{pro3.1}
Let $\beta >1$ and $0\leq \alpha<1$. The $\varphi$-expansion for
the dynamical system $\beta x+\alpha\mod1$ is valid.
$$
\Sigma^{\alpha,\beta}:=\big\{\i(x)\in\tA^{\Z_+}\colon x\in
X\backslash S\big\}=\big\{\ud{x}\in\tA^{\Z_+}\colon
\ud{u}^{\alpha,\beta}\prec\sigma^n\ud{x}\prec\ud{v}^{\alpha,\beta}\quad\forall
n\geq 0\big\}\,.
$$
Moreover
$$
\ud{u}^{\alpha,\beta}\preceq\sigma^n\ud{u}^{\alpha,\beta}\prec\ud{v}^{\alpha,\beta}
\quad\text{and}\quad
\ud{u}^{\alpha,\beta}\prec\sigma^n\ud{v}^{\alpha,\beta}\preceq\ud{v}^{\alpha,\beta}
\quad\forall n\geq0\,.
$$
\end{pro}

\noindent
The closure of $\Sigma^{\alpha,\beta}$ is the shift space
\begin{equation}\label{3.2}
\BSigma(\ud{u}^{\alpha,\beta},\ud{v}^{\alpha,\beta}):=\big\{\ud{x}\in\tA^{\Z_+}\colon
\ud{u}^{\alpha,\beta}\preceq\sigma^n\ud{x}\preceq\ud{v}^{\alpha,\beta}\quad\forall
n\geq 0\big\}\,.
\end{equation}
We define the {\sf orbits of $0$, resp. $1$} as, (the limits are
taken with $x\in X\backslash S$)
$$
T^k_\ab(0):=\lim_{x\downarrow 0}T_\ab^k(x)\,,\,k\geq
0\quad\text{resp.}\quad T^k_\ab(1):=\lim_{x\uparrow
1}T_\ab^k(x)\,,\,k\geq 0\,.
$$
From \eqref{2.12} and \eqref{2.12b} the coding of these orbits is
$\ud{u}^\ab$, resp. $\ud{v}^\ab$,
\begin{equation}\label{virtual}
\sigma^k\ud{u}^\ab=\lim_{x\downarrow
0}\i(T_\ab^k(x))\quad\text{and}\quad
\sigma^k\ud{v}^\ab=\lim_{x\uparrow 1}\i(T_\ab^k(x))\,.
\end{equation}
Notice that $T^k_\ab(0)<1$ and $T^k_\ab(1)>0$ for all $k\geq 0$.

The virtual itineraries $\ud{u}\equiv\ud{u}^\ab$ and
$\ud{v}\equiv\ud{v}^\ab$ of the dynamical system $\beta x+\alpha
\mod 1$  verify the conditions
\begin{equation}\label{condition}
\ud{u}\preceq\sigma^n\ud{u}\preceq \ud{v}\quad \forall\,n\geq
0\quad\text{and}\quad \ud{u}\preceq \sigma^n\ud{v}\preceq
\ud{v}\quad \forall\,n\geq 0\,.
\end{equation}
By Theorem \ref{thm2.1ter}, \eqref{virtual} and Theorem
\ref{thm2.5} we have $(x\in X\backslash S)$
\begin{align}
\ophi^\ab(\sigma^k\ud{u})&=\lim_{x\downarrow
0}\ophi^\ab(\i(T^k_\ab(x))=\lim_{x\downarrow 0}T^k_\ab(x)\equiv
T^k_\ab(0) \label{validity1}\\
\ophi^\ab(\sigma^k\ud{v})&=\lim_{x\uparrow
1}\ophi^\ab(\i(T^k_\ab(x))=\lim_{x\uparrow 1}T^k_\ab(x)\equiv
T^k_\ab(1)\,. \label{validity2}
\end{align}
Hence $\ud{u}$ and $\ud{v}$ verify the equations\footnote{If the
$\varphi$-expansion is not valid, which happens when $\beta=1$
and $\alpha\in\Q$, then \eqref{validity1} and \eqref{validity2}
are not necessarily true, as simple examples show. Hence
$\ud{u}^\ab$ and $\ud{v}^\ab$ do not necessarily verify
\eqref{equation}.}
\begin{equation}\label{equation}
\ophi^\ab(\ud{u})=0\,,\;
\ophi^\ab(\sigma\ud{u})=\alpha\quad\text{and}\quad
\ophi^\ab(\ud{v})=1\,,\; \ophi^\ab(\sigma\ud{v})=\gamma\,,
\end{equation}
with
\begin{equation}\label{gamma}
\gamma:=\alpha+\beta-k+1\in(0,1]\,.
\end{equation}
The strings $\ud{u}^\ab$ and $\ud{v}^\ab$ are
$\varphi$-expansions   of $0$ and $1$. Because of the presence of
discontinuities for the transformation $T_\ab$ at $a_1,\ldots
a_{k-1}$, there are other strings $\ud{u}$, $\ud{v}$ which verify
\eqref{condition} and \eqref{equation},  and which are also
$\varphi$-expansions   of $0$ and $1$. For latter purposes we
need to decribe these strings; this is the content of Proposition
\ref{pro3.4}, Proposition \ref{pro3.4bis} and Proposition
\ref{pro3.4ter}. We also take into consideration the borderline
cases $\alpha=1$ and $\gamma=0$. When $\alpha=1$ or $\gamma=0$
the dynamical system $T_\ab$ is defined using formula
\eqref{2.7}. The orbits of $0$ and $1$ are defined as before. For
example, if $\alpha=1$ it is the same dynamical system as
$T_{0,\beta}$, but with different symbols for the coding of the
orbits. The orbit of $0$ is coded by
$\ud{u}^{1,\beta}=(1)^\infty$, that is $\ud{u}^{1,\beta}_j=1$ for
all $j\geq 0$. Similarly, if $\gamma=0$ the orbit of $1$ is coded
by $\ud{v}^{\alpha,\beta}=(k-2)^\infty$. {\em We always assume
that $\alpha\in[0,1]$, $\gamma\in[0,1]$ and $\beta\geq1 $}.

\begin{lem}\label{lemelementary}
The equation
$$
y=\oph^\ab(x_k+t)\,, \quad y\in[0,1]
$$
can be solved uniquely if  $y\not\in S_0$, and its solution  is
$x_k=\i_0(y)$ and $t=T_\ab(y)\in(0,1)$. \\
If $y<y^\prime$, then the solutions of the equations
$$
y=\oph^\ab(x_k+t)\quad\text{and}\quad
y^\prime=\oph^\ab(x^\prime_k+t^\p)
$$
are such that either $x_k=x^\prime_k$ and
$T_\ab(y^\prime)-T_\ab(y)=\beta(y^\prime-y)$, or $x_k<x^\prime_k$.
\end{lem}

\prf The proof is elementary. It suffices to notice that
$$
y\not\in S_0\implies y= \varphi^\ab(x_k+t)\,.
$$
The second statement follows by monotonicity. \qed

\begin{pro}\label{pro3.4}
Let $0\leq\alpha< 1$ and assume that the $\varphi$-expansion is
valid. The following
assertions are equivalent.\\
1)\, There is a unique solution ($\ud{u}=\ud{u}^\ab$) of the
equations
\begin{equation}\label{eqalpha}
\ophi^\ab(\ud{u})=0\quad\text{and}\quad\ophi^\ab(\sigma\ud{u})=\alpha\,.
\end{equation}
2)\, The orbit of $0$ is not periodic or $x=0$ is a fixed point of $T_\ab$.\\
3)\, $\ud{u}^\ab$ is not periodic or $\ud{u}^\ab=\ud{0}$, where
$\ud{0}$ is the string $\ud{x}$ with $x_j=0$ $\forall j\geq 0$.
\end{pro}

\begin{pro}\label{pro3.4bis}
Let $0<\gamma\leq 1$ and assume that the $\varphi$-expansion is
valid. The following
assertions are equivalent.\\
1)\, There is a unique solution ($\ud{v}=\ud{v}^\ab$) of the
equations
\begin{equation}\label{eqbeta}
\ophi^\ab(\ud{v})=1\quad\text{and}\quad\ophi^\ab(\sigma\ud{v})=\gamma\,.
\end{equation}
2)\, The orbit of $1$ is not periodic or $x=1$ is a fixed point of $T_\ab$.\\
3)\, $\ud{v}^\ab$ is not periodic or $\ud{v}^\ab=(k-1)^\infty$.
\end{pro}

\prf We prove Proposition \ref{pro3.4}. Assume 1. The validity of
the $\varphi$-expansion implies that $\ud{u}^\ab$ is a solution of
\eqref{eqalpha}. If $\alpha=0$, then $\ud{u}^{0,\beta}=\ud{0}$ is
the only solution of \eqref{eqalpha} since  $\ud{x}\not=\ud{0}$
implies $\ophi^{0,\beta}(\ud{x})>0$ and $x=0$ is a fixed point of
$T_{0,\beta}$. Let $0<\alpha<1$.  Using Lemma \ref{lemelementary}
we deduce that $u_0=0$ and
$$
\alpha=T_\ab(0)=\oph^\ab(u_1+\ophi^\ab(\sigma^2\ud{u}))\,.
$$
If $\alpha=a_j$, $j=1,\ldots,k-1$ (see \eqref{eqso}), then
\eqref{eqalpha} has at least two solutions, which are
$0j(\sigma^2\ud{u}^\ab)$ with
$\ophi^\ab(\sigma^2\ud{u}^\ab)=T^2(0)=0$ (see \eqref{validity1}),
and $0(j-1)\ud{v}^\ab$ with $\ophi^\ab(\ud{v}^\ab)=1$. Therefore,
by our hypothesis we have $\alpha\not\in \{a_1,\ldots,a_{k-1}\}$,
$u_1=u_1^\ab$ and $\ophi^\ab(\sigma^2\ud{u}^\ab)=T^2(0)\in (0,1)$.
Iterating this argument we conclude that $1\implies 2$.\\
Assume 2. If $x=0$ is a fixed point, then $\alpha=0$ and
$\ud{u}^{0,\beta}=\ud{0}$. If the orbit of $0$ is not periodic,
\eqref{virtual} and the validity of the $\varphi$-expansion imply
$$
\sigma^k\ud{u}^\ab=\lim_{x\downarrow
0}\i(T^k_\ab(x))\succ\lim_{x\downarrow 0}\i(x)=\ud{u}^\ab\,.
$$
Assume 3. From \eqref{validity1} and the validity of the
$\varphi$-expansion we get
$$
\ophi^\ab(\sigma^k\ud{u}^\ab)=T^k_\ab(0)>
\ophi^\ab(\ud{u}^\ab)=0\,,
$$
so that the orbit of $0$ is not periodic.  The orbit of $0$ is not
periodic if and only if
$T^k_{\ab}(0)\not\in\{a_1,\ldots,a_{k-1}\}$ for all $k\geq 1$.
Using Lemma \ref{lemelementary} we conclude that \eqref{eqalpha}
has a unique solution. \qed

Propositions \ref{pro3.4} and \ref{pro3.4bis} give necessary and
sufficient conditions for the existence and uniqueness of the
solution of equations \eqref{equation}. In the following
discussion we  consider the case when there are several solutions.
The main results are summarize in Proposition \ref{pro3.4ter}. We
assume the validity of the $\varphi$-expansion.

Suppose first that the orbit of $1$ is not periodic and that the
orbit of $0$ is periodic, with minimal period $p:=\min\{k\colon
T^k(0)=0\}>1$. Hence $0<\gamma<1$ and $0<\alpha<1$. Let $\ud{u}$
be a solution of equations \eqref{eqalpha} and suppose furthermore
that $\ud{w}$ is a $\varphi$-expansion of $1$ such that
$$
\forall n\colon\; \ud{u}\preceq\sigma^n\ud{u}\preceq
\ud{w}\quad\text{with}\quad \ophi^\ab(\ud{w})=1\;,\;
\ophi^\ab(\sigma\ud{w})\leq\gamma\,.
$$
By Lemma \ref{lemelementary} we conclude that
$$
u_j=u_j^\ab\quad\text{and}\quad
T_\ab^{j+1}(0)=\ophi^\ab(\sigma^{j+1}\ud{u})\,,\quad
j=1,\ldots,p-2\,.
$$
Since $T^p(0)=0$, $T^{p-1}(0)\in \{a_1,\ldots,a_{k-1}\}$ and the
equation
$$
T_\ab^{p-1}(0)=\ophi^\ab\big(u_{p-1}+\ophi^\ab(\sigma^{p}\ud{u})\big)
$$
has two solutions.  Either $u_{p-1}=u_{p-1}^\ab$ and
$\ophi^\ab(\sigma^{p}\ud{u})=0$ or $u_{p-1}=u_{p-1}^\ab-1$ and
$\ophi^\ab(\sigma^{p}\ud{u})=1$. Let $\ud{a}$ be the prefix of
$\ud{u}^\ab$ of length $p$ and $\ud{a}^\prime$ the word of length
$p$ obtained by changing the last letter of $\ud{a}$
into\footnote{$u_{p-1}^\ab\geq 1$. $u_{p-1}^\ab=0$ if and only if
$p=1$ and $\alpha=0$.}
 $u_{p-1}^\ab-1$. We have
$\ud{a}^\prime<\ud{a}$. If $u_{p-1}=u_{p-1}^\ab$, then we can
again determine uniquely the next $p-1$ letters $u_i$. The
condition $\ud{u}\leq\sigma^k\ud{u}$ for  $k=p$ implies that we
have $u_{2p-1}=u_{p-1}^\ab$ so that, by iteration, we get the
solution $\ud{u}=\ud{u}^\ab$ for the equations \eqref{eqalpha}. If
$u_{p-1}=u_{p-1}^\ab-1$, then
$$
1=\ophi^\ab(\sigma^{p}\ud{u})=\ophi^\ab\big(u_p+\ophi^\ab(\sigma^{p+1}\ud{u})\big)\,.
$$
When $\ophi^\ab(\sigma^{p}\ud{u})=1$,  by our hypothesis on
$\ud{u}$ we also have $\ophi^\ab(\sigma^{p+1}\ud{u})=\gamma$. By
Proposition \ref{pro3.4bis} the equations
$$
\ophi^\ab(\sigma^{p}\ud{u})=1\quad\text{and}\quad
\ophi^\ab(\sigma^{p+1}\ud{u})=\gamma
$$
have  a unique solution, since we assume that the orbit of $1$ is
not periodic. The solution is $\sigma^{p}\ud{u}=\ud{v}^\ab$, so
that $\ud{u}=\ud{a}^\prime\ud{v}^\ab\prec\ud{u}^\ab$ is also a
solution of \eqref{eqalpha}. In that case there is no other
solution for \eqref{eqalpha}. The borderline case $\alpha=1$
corresponds to the periodic orbit of the fixed point $0$,
$\ud{u}^{1,\beta}=(1)^\infty$. Notice that
$\ophi^{1,\beta}(\sigma\ud{u}^{1,\beta})\not=1$. We can also
consider $\ophi^{1,\beta}$-expansions of $0$ with $u_0=0$ and
$\ophi^{1,\beta}(\sigma\ud{u})=1$. Our hypothesis on $\ud{u}$
imply that $\ophi^{1,\beta}(\sigma^2\ud{u})=\gamma$. Hence,
$\ud{u}=0\ud{v}^{1,\beta}=\ud{a}^\prime\ud{v}^\ab\prec\ud{u}^\ab$
is a solution of \eqref{eqalpha} and a
$\ophi^{1,\beta}$-expansion of $0$.

We can treat similarly the case when $\ud{u}^\ab$ is not periodic,
but $\ud{v}^\ab$ is periodic. When both $\ud{u}^\ab$ and
$\ud{v}^\ab$ are periodic we have more solutions, but the
discussion is similar. Assume that $\ud{u}^\ab$ has (minimal)
period $p>1$ and $\ud{v}^\ab$ has (minimal) period $q>1$. Define
$\ud{a}$, $\ud{a}^\prime$ as before, $\ud{b}$ as the prefix of
length $q$ of $\ud{v}^\ab$, and $\ud{b}^\prime$ as the word of
length $q$ obtained by changing the last letter of $\ud{b}$ into
$v_{q-1}^\ab+1$. When $0<\alpha<1$ and $0<\gamma<1$, one shows as
above that the elements $\ud{u}\not=\ud{u}^\ab$ and
$\ud{v}\not=\ud{v}^\ab$ which are $\oph^\ab$-expansions of $0$ and
$1$ are of the form
$$
\ud{u}=\ud{a}^\prime\ud{b}^{n_1}\ud{b}^\prime\ud{a}^{n_2}\cdots\,,\,n_i\geq
0 \quad\text{and}\quad
\ud{v}=\ud{b}^\prime\ud{a}^{m_1}\ud{a}^\prime\ud{b}^{m_2}\cdots\,,\,m_i\geq
0 \,.
$$
The integers $n_i$ and $m_i$ must be such that \eqref{condition}
is verified. The largest solution of \eqref{eqalpha} is
$\ud{u}^\ab$ and the smallest one  is $\ud{a}^\prime\ud{v}^\ab$.

\begin{pro}\label{pro3.4ter}
Assume that  the $\varphi$-expansion is valid.\\
1)\, Let $\ud{u}$ be a solution of \eqref{eqalpha}, such that
$\ud{u}\preceq\sigma^n\ud{u}$ for all $n\geq 1$, and let $\ud{v}$
be a solution of \eqref{eqbeta}, such that
$\sigma^n\ud{v}\preceq\ud{v}$ for all $n\geq 1$. Then
$$
\ud{u}\preceq\ud{u}^\ab\quad\text{and}\quad
\ud{v}^\ab\preceq\ud{v}\,.
$$
 2)\, Let  $\ud{u}$ be a solution of
\eqref{eqalpha}, and let $\ud{u}^\ab=(\ud{a})^\infty$ be periodic
with minimal period $p>1$, and suppose that there exists $\ud{w}$
such that
$$
\forall n\colon\; \ud{u}\preceq\sigma^n\ud{u}\preceq
\ud{w}\quad\text{with}\quad \ophi^\ab(\ud{w})=1\;,\;
\ophi^\ab(\sigma\ud{w})\leq\gamma\,.
$$
Then
\begin{equation}\label{3.3.3}
\ud{u}^{\ab}_*\preceq\ud{u}\preceq\ud{u}^\ab\quad\text{where}\quad
\ud{u}^{\ab}_*:=\ud{a}^\prime\ud{v}^\ab\;\text{and}\;
\ud{a}^\p:=(\tp\ud{a})(a_{p-1}-1)\,.
\end{equation}
Moreover, $\ud{u}=\ud{u}^\ab$ $\iff$ $\ud{a}$ is a prefix
of $\ud{u}$ $\iff$ $\ophi^\ab(\sigma^p\ud{u})<1$. \\
3)\, Let $\ud{v}$ be a solution of \eqref{eqbeta}, and let
$\ud{v}^\ab=(\ud{b})^\infty$ be periodic with minimal period
$q>1$, and suppose that there exists $\ud{w}$ such that
$$
\forall n\colon\; \ud{w}\preceq\sigma^n\ud{v}\preceq
\ud{v}\quad\text{with}\quad \ophi^\ab(\ud{w})=0\;,\;
\ophi^\ab(\sigma\ud{w})\geq\alpha\,.
$$
Then
\begin{equation}\label{3.3.4}
\ud{v}^\ab\preceq\ud{v}\preceq
\ud{v}^{\ab}_*\quad\text{where}\quad
\ud{v}^{\ab}_*:=\ud{b}^\prime\ud{u}^\ab\;\text{and}\;
\ud{b}^\p:=(\tp\ud{b})(b_{q-1}+1)\,.
\end{equation}
Moreover, $\ud{v}=\ud{v}^\ab$ $\iff$ $\ud{b}$ is a prefix of
$\ud{v}$ $\iff$ $\ophi^\ab(\sigma^q\ud{v})>0$.
\end{pro}

\section{Shift space ${\mathbf\Sigma}(\ud{u},\ud{v})$}\label{section3}
\setcounter{equation}{0}

Let $\ud{u}\in\tA^{\Z_+}$ and $\ud{v}\in\tA^{\Z_+}$, such that
$u_0=0$, $v_0=k-1$ ($k\geq 2$) and \eqref{condition} holds. These
assumptions are valid for the whole section, except subsection
\ref{subsectionalgo}. We study the shift-space
\begin{equation}\label{3.2.2}
\BSigma(\ud{u},\ud{v}):=\{\ud{x}\in\tA^{\Z_+}\colon \ud{u}\preceq
\sigma^n\ud{x}\preceq \ud{v}\quad \forall\,n\geq 0\}\,.
\end{equation}
It is useful to extend the relation $\prec$ to words or to words
and strings. We do it only in the following case. Let $\ud{a}$ and
$\ud{b}$ be words (or strings). Then
$$
\ud{a}\prec\ud{b}\quad\text{iff $\exists$ $\ud{c}\in\tA^*$,
$\exists$ $k\geq 0$ such that $\ud{a}=\ud{c}a_k\cdots$,
$\ud{b}=\ud{c}b_k\cdots$ and $a_k<b_k$.}
$$
If $\ud{a}\prec\ud{b}$ then neither $\ud{a}$ is a prefix of
$\ud{b}$, nor $\ud{b}$ is a prefix of $\ud{a}$.

In subsection \ref{subsectionfollower} we introduce one of the
main tool for studying the shift-space $\BSigma(\ud{u},\ud{v})$,
the follower-set graph.  In subsection \ref{subsectionalgo} we
give an algorithm which assigns to a pair of strings
$(\ud{u},\ud{v})$ a pair of real numbers
$(\bar{\alpha},\bar{\beta})\in [0,1]\times[1,\infty)$. Finally in
subsection \ref{topological} we compute the topological entropy
of the shift space $(\ud{u},\ud{v})$.

\subsection{Follower-set graph
$\cG(\ud{u},\ud{v})$}\label{subsectionfollower}

We associate to $\BSigma(\ud{u},\ud{v})$ a graph
$\cG(\ud{u},\ud{v})$, called the {\sf follower-set graph} (see
\cite{LiM}), as well as an equivalent graph
$\overline{\cG}(\ud{u},\ud{v})$. The graph
$\overline{\cG}(\ud{u},\ud{v})$  has been systematically studied
by Hofbauer in his works about piecewise monotone one-dimensional
dynamical systems; see  \cite{Ho1}, \cite{Ho2} and \cite{Ho3} in
the context of this paper, as well as  \cite{Ke} and \cite{BrBr}.
Our presentation differs from that of Hofbauer, but several proofs
are directly inspired by \cite{Ho2} and \cite{Ho3}.

We denote by $\cL(\ud{u},\ud{v})$ the {\sf language of}
$\BSigma(\ud{u},\ud{v})$, that is the set of words, which are
factors of $\ud{x}\in \BSigma(\ud{u},\ud{v})$ (including the empty
word $\epsilon$). Since $\sigma\BSigma(\ud{u},\ud{v})\subset
\BSigma(\ud{u},\ud{v})$, the language is also the set of prefixes
of the strings $\ud{x}\in \BSigma(\ud{u},\ud{v})$.  To simplify
the notations we set in this subsection
$\BSigma:=\BSigma(\ud{u},\ud{v})$, $\cL:=\cL(\ud{u},\ud{v})$,
$\cG:=\cG(\ud{u},\ud{v})$.

Let $\cC_u$ be the set of words $\ud{w}\in\cL$ such that
$$
\ud{w}=\begin{cases} \text{$\ud{w}^\prime$\colon
$\ud{w}^\prime\not=\epsilon$, $\ud{w}^\prime$ is a prefix of
$\ud{u}$}\\
\text{$w_0\ud{w}^\prime$\colon $w_0\not=u_0$, $\ud{w}^\prime$ is a
prefix of $\ud{u}$, possibly $\epsilon$.}
\end{cases}
$$
Similarly we introduce $\cC_v$ as the set of words $\ud{w}\in\cL$
such that
$$
\ud{w}=\begin{cases} \text{$\ud{w}^\prime$\colon
$\ud{w}^\prime\not=\epsilon$, $\ud{w}^\prime$ is a prefix of
$\ud{v}$}\\
\text{$w_0\ud{w}^\prime$\colon $w_0\not=v_0$, $\ud{w}^\prime$ is a
prefix of $\ud{v}$, possibly $\epsilon$.}
\end{cases}
$$

\begin{defn}\label{defn3.1}
Let $\ud{w}\in\cL$. The longest suffix of $\ud{w}$, which is a
prefix of $\ud{v}$, is denoted by $v(\ud{w})$. The longest suffix
of $\ud{w}$, which is a prefix of $\ud{u}$, is denoted by
$u(\ud{w})$. The {\sf $u$-parsing of $\ud{w}$} is the following
decomposition of $\ud{w}$ into $\ud{w}=\ud{a}^1\cdots\ud{a}^k$
with $\ud{a}^j\in\cC_u$. The first word $\ud{a}^1$ is the longest
prefix of $\ud{w}$ belonging to $\cC_u$. If
$\ud{w}=\ud{a}^1\ud{w}^\prime$ and $\ud{w}^\prime\not=\epsilon$,
then the next word $\ud{a}^2$ is the longest prefix of
$\ud{w}^\prime$ belonging to $\cC_u$ and so on.

The {\sf $v$-parsing of $\ud{w}$} is the analogous decomposition
of $\ud{w}$ into $\ud{w}=\ud{b}^1\cdots\ud{b}^\ell$ with
$\ud{b}^j\in\cC_v$.
\end{defn}

\begin{lem}\label{lem3.2.0}
Let $\ud{w}\ud{c}$ and $\ud{c}\ud{w}^\prime$ be prefixes of
$\ud{u}$ (respectively of $\ud{v}$). If
$\ud{w}\ud{c}\ud{w}^\prime\in\cL$, then
$\ud{w}\ud{c}\ud{w}^\prime$ is a prefix of $\ud{u}$ (respectively
of $\ud{v}$). Let $\ud{w}\in\cL$. If $\ud{a}^1\cdots\ud{a}^k$ is
the $u$-parsing of $\ud{w}$, then only the first word $\ud{a}^1$
can be a prefix of $\ud{u}$, otherwise $u(\ud{a}^j)=\s\ud{a}^j$.
Moreover $u(\ud{a}^k)=u(\ud{w})$. Analogous properties hold for
the $v$-parsing of $\ud{w}$.
\end{lem}

\prf Suppose that $\ud{w}\ud{c}$ and $\ud{c}\ud{w}^\prime$ are
prefixes of $\ud{u}$. Then $\ud{w}$ is a prefix of $\ud{u}$.
Assume that $\ud{w}\ud{c}\ud{w}^\prime\in\cL$ is not a prefix of
$\ud{u}$. Then $\ud{u}\prec\ud{w}\ud{c}\ud{w}^\prime$. Since
$\ud{w}$ is a prefix of $\ud{u}$,
$\sigma^{|\ud{w}|}\ud{u}\prec\ud{c}\ud{w}^\prime$. This
contradicts the fact that $\ud{c}\ud{w}^\prime$ is a prefix of
$\ud{u}$. By applying this result with $\ud{c}=\epsilon$ we get
the result that only the first word in the $u$-parsing of
$\ud{w}$ can be a prefix of $\ud{u}$. Suppose that the
$u$-parsing of $\ud{w}$ is $\ud{a}^1\cdots\ud{a}^k$. Let $k\geq
2$ and assume that $u(\ud{w})$ is not a suffix of $\ud{a}^k$ (the
case $k=1$ is obvious). Since $\ud{a}^k$ is not a prefix of
$\ud{u}$, $u(\ud{w})$ has $\ud{a}^k$ as a proper suffix. By the
first part of the lemma this contradicts the maximality property
of the words in the $u$-parsing. \qed

\begin{lem}\label{lem3.2.1}
Let $\ud{w}\in\cL$. Let $p=|u(\ud{w})|$ and $q=|v(\ud{w})|$. Then
$$
\big\{\ud{x}\in\BSigma\colon \text{$\ud{w}$ is a prefix of
$\ud{x}$}\big\}= \big\{\ud{x}\in \tA^{\Z_+}\colon
\ud{x}=\ud{w}\ud{y}\,,\;\ud{y}\in\BSigma\,,\;
\sigma^p\ud{u}\preceq\ud{y}\preceq\sigma^q\ud{v}\big\}\,.
$$
Moreover,
$$
\big\{\ud{y}\in\BSigma\colon
\ud{w}\ud{y}\in\BSigma\big\}=\big\{\ud{y}\in\BSigma\colon
u(\ud{w})\ud{y}\in\BSigma\big\}\quad\text{if $p>q$}
$$
$$
\;\big\{\ud{y}\in\BSigma\colon
\ud{w}\ud{y}\in\BSigma\big\}=\big\{\ud{y}\in\BSigma\colon
v(\ud{w})\ud{y}\in\BSigma\big\}\quad\text{if $q>p$.}
$$
\end{lem}

\prf Suppose that $\ud{x}\in\BSigma$ and $\ud{w}$, $|\ud{w}|=n$,
is a prefix of $\ud{x}$.  Let $n\geq 1$ (the case $n=0$ is
trivial). We can write $\ud{x}=\ud{w}\ud{y}$. Since
$\ud{x}\in\BSigma$,
$$
\ud{u}\preceq
\sigma^{\ell+n}\ud{x}\preceq\ud{v}\quad\forall\,\ell\geq 0\,,
$$
so that $\ud{y}\in\BSigma$. We  have
$$
\ud{u}\preceq \sigma^{n-p}\ud{x}= u(\ud{w})\ud{y}\,.
$$
Since $u(\ud{w})$ is a prefix of $\ud{u}$ of length $p$, we get
$\sigma^p\ud{u}\preceq \ud{y}$. Similarly we prove that
$\ud{y}\preceq\sigma^q\ud{v}$.

Suppose that $\ud{x}=\ud{w}\ud{y}$, $\ud{y}\in\BSigma$ and
$\sigma^p\ud{u}\preceq\ud{y}\preceq\sigma^q\ud{v}$. To prove that
$\ud{x}\in\BSigma$, it is sufficient to prove that
$\ud{u}\preceq\sigma^m\ud{x}\preceq \ud{v}$ for $m=0,\ldots, n-1$.
We prove $\ud{u}\preceq\sigma^m\ud{x}$ for $m=0,\ldots, n-1$. The
other case is similar. Let $\ud{w}=\ud{a}^1\cdots\ud{a}^\ell$ be
the $u$-parsing of $\ud{w}$, $|\ud{w}|=n$  and $p=|u(\ud{w})|$. We
have
$$
\sigma^p\ud{u}\preceq\ud{y}\implies \ud{u}\preceq
\sigma^j\ud{u}\preceq \sigma^ju(\ud{w})\ud{y}\quad\forall
\,j=0,\ldots,p\,.
$$
If $\ud{a}^\ell$ is not a prefix of $\ud{u}$, then  $p=n-1$ and
we also have  $\ud{u}\preceq \ud{a}^k\ud{y}$. If $\ud{a}^\ell$ is
a prefix of $\ud{u}$, then  $p=n$ (and $\ell=1$). This proves the
result for $\ell=1$. Let $\ell\geq 2$. Then $\ud{a}^\ell$ is not
a prefix of $\ud{u}$ and $\ud{a}^{\ell-1}\ud{a}^\ell\in\cL$.
Suppose that $\ud{a}^{\ell-1}$ is not a prefix of $\ud{u}$. In
that case $\ud{u}\preceq\ud{a}^{\ell-1}\ud{a}^\ell\ud{y}$ and we
want to prove that
$\ud{u}\preceq\sigma^{j}\ud{a}^{\ell-1}\ud{a}^\ell\ud{y}$ for
$j=1,\ldots, |\ud{a}^{\ell-1}|$. We know that
$\sigma\ud{a}^{\ell-1}$ is a prefix of $\ud{u}$, and by maximality
of the words in the $u$-parsing and Lemma \ref{lem3.2.0}
$\ud{u}\prec\sigma\ud{a}^{\ell-1}\ud{a}^\ell$; hence
$\ud{u}\prec\sigma\ud{a}^{\ell-1}\ud{a}^\ell\ud{y}$. Therefore
$$
\ud{u}\preceq \sigma^j\ud{u}\preceq
\sigma^j\ud{a}^{\ell-1}\ud{a}^\ell\ud{y}\quad\forall
\,j=0,\ldots,|\ud{a}^{\ell-1}|\,.
$$
Similar proof if $\ell=2$ and $\ud{a}^{\ell-1}$ is a prefix of
$\ud{u}$. Iterating this argument we prove that
$\ud{u}\preceq\sigma^m\ud{x}$ for $m=0,\ldots,n-1$.

Suppose that $|u(\ud{w})|>|v(\ud{w})|$ and set $\ud{a}=u(\ud{w})$.
We prove that $v(\ud{a})=v(\ud{w})$. By definition $v(\ud{w})$ is
the longest suffix of $\ud{w}$ which is a prefix of $\ud{v}$; it
is also a suffix of $\ud{a}$, whence it is also the longest suffix
of $\ud{a}$ which is a prefix of $\ud{v}$. Therefore, from the
first part of the lemma we get
$$
\big\{\ud{y}\in\BSigma\colon
\ud{w}\ud{y}\in\BSigma\big\}=\big\{\ud{y}\in\BSigma\colon
u(\ud{w})\ud{y}\in\BSigma\big\}\,.
$$
 \qed

\begin{defn}\label{defn3.2.1}
Let $\ud{w}\in\cL$. The {\sf follower-set}\footnote{Usually the
follower-set is defined as $\cF_{\ud{w}}=\big\{\ud{y}\in\cL\colon
\ud{w}\ud{y}\in\cL\big\}$. Since $\cL$ is a dynamical language,
i.e. for each $\ud{w}\in\cL$ there exists a letter $e\in\tA$ such
that $\ud{w}e\in\cL$, the two definitions agree.} of $\ud{w}$ is
the set
$$
\cF_{\ud{w}}:=\big\{\ud{y}\in\BSigma\colon
\ud{w}\ud{y}\in\BSigma\big\}\,.
$$
\end{defn}

Lemma \ref{lem3.2.1} gives the important results that
$\cF_{\ud{w}}=\cF_{u(\ud{w})}$ if $|u(\ud{w})|>|v(\ud{w})|$, and
$\cF_{\ud{w}}=\cF_{v(\ud{w})}$ if $|v(\ud{w})|>|u(\ud{w})|$.
Moreover,
\begin{equation}\label{3.2.3}
\cF_\ud{w}=\big\{\ud{y}\in\BSigma\colon
\sigma^p\ud{u}\preceq\ud{y}\preceq\sigma^q\ud{v}\big\}\quad\text{where
$p=|u(\ud{w})|$ and $q=|v(\ud{w})|$.}
\end{equation}
We can define an equivalence relation between words of $\cL$,
$$
\ud{w}\sim \ud{w}^\prime \iff \cF_{\ud{w}}=\cF_{\ud{w}^\prime}\,.
$$
The collection of follower-sets is entirely determined by the
strings $\ud{u}$ and $\ud{v}$. Moreover, the strings $\ud{u}$ and
$\ud{v}$ are eventually periodic if and only if this collection
is finite. Notice that $\BSigma=\cF_\epsilon=\cF_{\ud{w}}$ when
$p=q=0$.

\begin{defn}\label{defn3.2.2}
The {\sf follower-set graph} $\cG$ is the labeled graph whose set
of vertices is the collection of all follower-sets. Let $\cC$ and
$\cC^\prime$ be two vertices. There is an edge, labeled by
$a\in\tA$, from $\cC$ to $\cC^\prime$ if and only if there exists
$\ud{w}\in\cL$ so that $\ud{w}a\in\cL$, $\cC=\cF_{\ud{w}}$ and
$\cC^\prime=\cF_{\ud{w}a}$. $\cF_\epsilon$ is called the {\sf
root} of $\cG$.
\end{defn}

The following properties of $\cG$ are immediate. From any vertex
there is at least one out-going edge and at most $|\tA|$. If
$\tA=\{0,1,\ldots,k-1\}$ and $k\geq 3$, then for each $j\in
\{1,\ldots,k-2\}$ there is an edge labeled by $j$ from
$\cF_\epsilon$ to $\cF_\epsilon$. The out-going edges from
$\cF_{\ud{w}}$ are labeled by the first letters of the strings
$\ud{y}\in \cF_{\ud{w}}$. The follower-set graph $\cG$ is
right-resolving. Given $\ud{w}\in\cL$, there is a unique path
labeled by $\ud{w}$  from  $\cF_\epsilon$ to $\cF_{\ud{w}}$.

\begin{lem}\label{lem3.2.2}
Let $\ud{a}$ be a $u$-prefix and suppose that $\ud{b}=v(\ud{a})$.
Let $p=|\ud{a}|$ and $q=|\ud{b}|$ so that
$\cF_\ud{a}=\big\{\ud{y}\in\BSigma\colon
\sigma^p\ud{u}\preceq\ud{y}\preceq\sigma^q\ud{v}\big\}$. Then
there are more than one out-going edges from $\cF_{\ud{a}}$ if
and only if $u_p<v_q$.

Assume that $u_p<v_q$. Then there is an edge labeled by $v_q$ from
$\cF_{\ud{a}}$ to $\cF_{\ud{b}v_q}$, an edge labeled by $u_p$
from $\cF_{\ud{a}}$ to $\cF_{\ud{a}u_p}$ and
$v(\ud{a}u_p\ud{c})=v(\ud{c})$. If there exists $u_p<\ell<v_q$,
there is an edge labeled by $\ell$ from $\cF_{\ud{a}}$ to
$\cF_\epsilon$. Moreover, there are at least two out-going edges
from $\cF_{\ud{b}}$, one labeled by $v_q$ to $\cF_{\ud{b}v_q}$
and one labeled by $\ell^\prime=u_{|u(\ud{b})|+1}<v_q$ to
$\cF_{u(\ud{b})\ell^{\prime}}$. Furthermore
$u(\ud{b}v_q\ud{c})=u(\ud{c})$.
\end{lem}

\prf The first part of the lemma is immediate. Suppose that there
is only one out-going edge from $\cF_{\ud{b}}$, that is from
$\cF_{\ud{b}}$ to $\cF_{\ud{b}v_q}$. This happens if and only if
$u(\ud{b})v_q$ is a prefix of $\ud{u}$. By Lemma \ref{lem3.2.0} we
conclude that $\ud{a}v_q$ is a prefix of $\ud{u}$, which is a
contradiction. Therefore $u(\ud{b}v_q)=\epsilon$; hence
$u(\ud{b}v_q\ud{c})=u(\ud{c})$. \qed

\begin{lem}\label{lem3.2.3}
Let $\ud{b}$ be a $v$-prefix and suppose that $\ud{a}=u(\ud{b})$.
Let $p=|\ud{a}|$ and $q=|\ud{b}|$ so that
$\cF_\ud{b}=\big\{\ud{y}\in\BSigma\colon
\sigma^p\ud{u}\preceq\ud{y}\preceq\sigma^q\ud{v}\big\}$. Then
there are more than one out-going edges from $\cF_{\ud{b}}$ if and
only if $u_p<v_q$.

Assume that $u_p<v_q$. Then there is an edge labeled by $u_p$ from
$\cF_{\ud{b}}$ to $\cF_{\ud{a}u_p}$, an edge labeled by $v_q$
from $\cF_{\ud{b}}$ to $\cF_{\ud{b}v_q}$ and
$u(\ud{b}v_q\ud{c})=u(\ud{c})$. If there exists $u_p<\ell<v_q$,
there is an edge labeled by $\ell$ from $\cF_{\ud{b}}$ to
$\cF_\epsilon$. Moreover, there are at least two out-going edges
from $\cF_{\ud{a}}$, one labeled by $u_p$ to $\cF_{\ud{a}u_p}$
and one labeled by $\ell^\prime=v_{|v(\ud{a})|+1}>u_p$ to
$\cF_{v(\ud{a})\ell^{\prime}}$. Furthermore
$v(\ud{a}u_p\ud{c})=v(\ud{c})$.
\end{lem}

\begin{sch}\label{sch}
 The picture below illustrates the main properties of the graph
$\cG$. The vertices of the graph are labeled by prefixes of
$\ud{u}$ and $\ud{v}$. The above line represents a prefix of the
string $\ud{u}$ which is written $\ud{w}\,e^{\p\p}$, and the
bottom line a prefix of the string $\ud{v}$, which is written
$\ud{b}\,e^\p$. Here $u(\ud{b})=\ud{a}$, $v(\ud{w})=\ud{b}$ and we
assume that $e^{\p\p}\prec e^\prime$. Therefore, there is an edge
labeled by $e^\prime$ from $\cF_{\ud{w}}$ to
$\cF_{\ud{b}e^\prime}$ and there is an edge labeled by $e$ from
$\cF_{\ud{b}}$ to $\cF_{\ud{a}e}$ with $e\prec e^\p$. Moreover,
we also have $e\preceq e^{\p\p}$. Only these two labeled edges are
drawn in the picture.

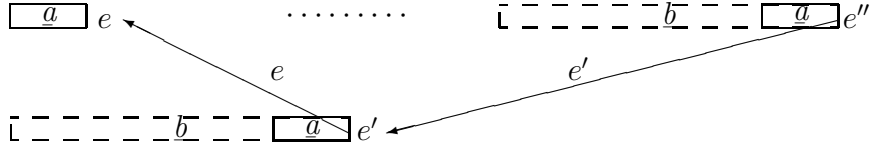
\setlength{\unitlength}{1mm}
\begin{picture}(160,35)
\put(15,25){\framebox(10,3)[b]{$\ud{a}$}}
\put(40,25){\makebox(40,3){$\cdots\cdots\cdots$}}
\put(80,25){\dashbox{2}(45,3){$\ud{b}$}}\put(115,25){\framebox(10,3){$\ud{a}$}}
\put(15,10){\dashbox{2}(45,3){$\ud{b}$}}\put(50,10){\framebox(10,3)[b]{$\ud{a}$}}
\put(60,11){\vector(-2,1){30}} \put(125,26){\vector(-4,-1){60}}
\put(25,25){\makebox[5mm]{$e$}} \put(60,10){\makebox[5mm]{$e^\p$}}
\put(125,25){\makebox[5mm]{$e^{\p\p}$}}
\put(48,18){\makebox[5mm]{$e$}} \put(88,18){\makebox[5mm]{$e^\p$}}
\end{picture}
\end{sch}

We introduce a variant of the follower-set graph denoted below by
$\overline{\cG}(\ud{u},\ud{v})$ or simply by $\overline{\cG}$. We
introduce a vertex for each (nontrivial) prefix $\ud{a}$ of
$\ud{u}$ and for each (nontrivial) prefix of $\ud{b}$ of
$\ud{v}$. We add the vertex $\cF_\epsilon$. Here we do not use
the equivalence relation $\sim$. The root is denoted by $[0,0]$;
let $\ud{a}$ be a prefix of $\ud{u}$ and let $p=|\ud{a}|$,
$q=|v(\ud{a})|$. Then this vertex is denoted by $[p,q]$. Notice
that $p>q$. Similarly, let $\ud{b}$ be a prefix of $\ud{v}$ and
let $p=|u(\ud{b})|$, $q=|\ud{b}|$. Then the corresponding vertex
is denoted by $[p,q]$. Notice that $p<q$. The {\sf upper branch}
of $\overline{\cG}$ is the set of all vertices $[p,q]$ with $p>q$
and the {\sf lower branch} of $\overline{\cG}$ is the set of all
vertices $[p,q]$ with $p<q$. There is a single out-going edge
from $[p,q]$ if and only if $u_p=v_q$. In that case the edge is
labeled by $u_p$ (or $v_q$) and goes from $[p,q]$ to $[p+1,q+1]$.
Otherwise there are several out-going edges. If  $u_p<v_q$ there
is an edge labeled by $u_p$ from $[p,q]$ to $[p+1,0]$, an edge
labeled by $v_q$ from $[p,q]$ to $[0,q+1]$, and if $u_p<j<v_q$
then there is an edge labeled by $j$ from $[p,q]$ to $[0,0]$. We
define the {\sf level of the vertex} $[p,q]$ of $\overline{\cG}$
as $\ell([p,q]):=\max\{p,q\}$.

\begin{defn}\label{defn3.2.3}
Let $\BSigma$ be a shift-space and $\cL$ its language. We denote
by $\cL_n$ the set of all words of $\cL$ of length $n$. The {\sf
 entropy} of $\BSigma$ is
$$
h(\BSigma):=\lim_{n\ra\infty}\frac{1}{n}\log_2 {\rm
card}(\cL_n)\,.
$$
\end{defn}

The number $h(\BSigma)$ is also equal to the topological entropy
of the dynamical system $(\BSigma,\sigma)$ \cite{LiM}. In our case
we can give an equivalent definition using the graph $\cG$ or the
graph $\overline{\cG}$.  We set
$$
\ell(n):={\rm card}\{\text{$n$-paths in $\overline{\cG}$ starting
at the root $\cF_\epsilon$}\}\,.
$$
Since the graph is right-resolving and for any $\ud{w}\in\cL_n$
there is a unique  path labeled by $\ud{w}$, starting at the root
$[0,0]$, so that $h(\BSigma)=h(\overline{\cG})$ where
$$
h(\overline{\cG})=\lim_{n\ra\infty}\frac{1}{n}\log_2 \ell(n)\,.
$$

Let $K\in\N$ and $\overline{\cG}_K$ be the sub-graph of
$\overline{\cG}$ whose set of vertices is the set of all vertices
of $\overline{\cG}$ of levels smaller or equal to $K$. The
following result is Proposition 9.3.15 in \cite{BrBr}.

\begin{pro}\label{proBB}
Given $\varepsilon>0$ there exists a $K(\varepsilon)<\infty$ such
that for any $K\geq K(\varepsilon)$,
$$
h(\overline{\cG}_K)\leq h(\overline{\cG})\leq
h(\overline{\cG}_K)+\varepsilon\,.
$$
\end{pro}

\begin{cor}\label{corBB}
Let $(\ud{u},\ud{v})$ be a pair of strings of $\tA^{\Z_+}$
verifying \eqref{condition}. Given $\varepsilon>0$ there exists
$N(\varepsilon)$ such that if $(\ud{u}^\prime,\ud{v}^\prime)$ is
a pair of strings verifying \eqref{condition}, $\ud{u}$ and
$\ud{u}^\prime$ have a common prefix of length larger
$N(\varepsilon)$ and $\ud{v}$ and $\ud{v}^\prime$ have a common
prefix of length larger than $N(\varepsilon)$, then
$$
\big|h(\BSigma(\ud{u}^\prime,\ud{v}^\prime))-
h(\BSigma(\ud{u},\ud{v}))\big|\leq\varepsilon\,.
$$
\end{cor}

\subsection{The algorithm for finding $(\bar{\alpha},\bar{\beta})$}\label{subsectionalgo}

We describe an algorithm, which assigns to a pair of strings
$(\ud{u},\ud{v})$, such that $u_0=0$ and $v_0=k-1$,  a pair of
real numbers $(\bar{\alpha},\bar{\beta})\in
[0,1]\times[1,\infty)$. We assume tacitly that for the pair
$(\ab)$ one has $\alpha\in[0,1]$, $\beta\leq k$, and that the map
$\oph^\ab$ verifies
$$
0<\oph^\ab(t)<1\quad\forall t\in(1,k-1)\,.
$$
In particular $\beta\geq k-2$. When $k=2$ we assume that
$\beta\geq 1$. Recall that
$$
\gamma=\alpha+\beta-k+1\,,
$$
and notice that our assumptions imply that $0\leq\gamma\leq 1$.

\begin{defn}\label{defndominate}
The map $\oph^\ab$ {\sf dominates} the map $\oph^\abp$ if and only
if $\oph^\ab(t)\geq \oph^\abp(t)$ for all $t\in [0,k]$ and there
exists $s\in [0,k]$ such that $\oph^\ab(s)>\oph^\abp(s)$.
\end{defn}

\setlength{\unitlength}{1mm}
\begin{picture}(160,50)
\put(25,15){\framebox(30,30)} \put(55,15){\framebox(30,30)}
\put(85,15){\framebox(30,30)} \put(35,15){\line(5,2){75}}
\put(25,23){\dashbox{2}(30,0)} \put(25,35){\dashbox{2}(60,0)}
\put(20,22){\makebox[5mm]{$a_1$}}
\put(20,34){\makebox[5mm]{$a_2$}}
\put(32,10){\makebox[5mm]{$\alpha$}}
\put(108,10){\makebox[5mm]{$2+\gamma=\alpha+\beta$}}
\put(34,14){$\bullet$}
\put(109,14){$\bullet$}
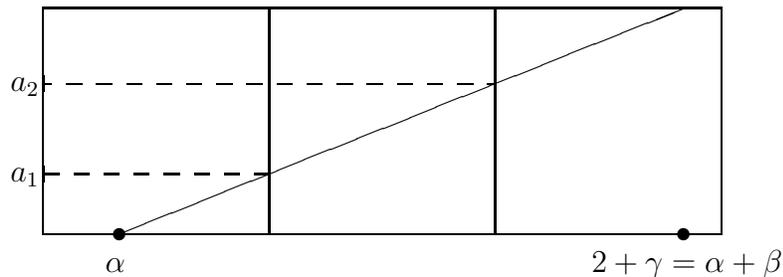
\put(55,5){\makebox[30mm]{\it{The graph of $\oph^\ab$} ($k=3$)}}
\end{picture}

\begin{lem}\label{lem3.3.1}
If $\oph^\ab$ dominates $\oph^\abp$, then, for all
$\ud{x}\in\tA^{\Z_+}$, $\ophi^\ab(\ud{x})\geq
\ophi^\abp(\ud{x})$. If
$$
0<\ophi^\ab(\ud{x})<1\quad\text{or}\quad 0<\ophi^\abp(\ud{x})<1\,,
$$
then the inequality is strict.
\end{lem}

\prf If $\oph^\ab$ dominates $\oph^\abp$, then by our implicit
assumptions we get by inspection of the graphs that
$$
\forall t\geq t^\p \colon\;\oph^\ab(t)>\oph^\abp(t^\p)\quad
\text{if}\quad t,t^\p\in (\alpha,\alpha^\prime+\beta^\prime)=
(\alpha,\alpha+\beta)\cup(\alpha^\prime,\alpha^\prime+\beta^\prime)\,,
$$
otherwise $\oph^\ab(t)\geq\oph^\abp(t^\p)$. Therefore, for all
$n\geq 1$,
$$
\oph_n^\ab(\ud{x})\geq\oph_n^\abp(\ud{x})\,.
$$
Suppose that $0<\ophi^\ab(\ud{x})<1$. Then
$x_0+\ophi^\ab(\sigma\ud{x})\in(\alpha,\alpha+\beta)$ and
$$
\ophi^\ab(\ud{x})=\oph^\ab(x_0+\ophi^\ab(\sigma\ud{x}))>
\oph^\abp(x_0+\ophi^\abp(\sigma\ud{x}))=\ophi^\abp(\ud{x})\,.
$$
Similar proof for $0<\ophi^\abp(\ud{x})<1$. \qed

\begin{lem}\label{lem3.3.2}
Let $\alpha=\alpha^\prime\in [0,1]$ and $1\leq
\beta<\beta^\prime$. Then, for
$\ud{x}\in\tA^{\Z_+}$,
$$
0\leq
\ophi^{\alpha,\beta}(\ud{x})-\ophi^{\alpha,\beta^\prime}(\ud{x})\leq
\frac{|\beta-\beta^\prime|}{\beta^\prime-1}\,.
$$
Let $\gamma=\gamma^\prime\in [0,1]$,
$0\leq\alpha^\prime<\alpha\leq 1$ and $\beta^\prime>1$. Then, for
$\ud{x}\in\tA^{\Z_+}$,
$$
0\leq \ophi^{\abp}(\ud{x})-\ophi^{\alpha,\beta}(\ud{x})\leq
\frac{|\alpha-\alpha^\prime|}{\beta^\prime-1}\,.
$$
The map $\beta\mapsto \ophi^{\alpha,\beta}(\ud{x})$ is continuous
at $\beta=1$.
\end{lem}

\prf Let $\alpha=\alpha^\prime\in [0,1]$ and $1\leq
\beta<\beta^\prime$.
For $t,t^\prime\in [0,k]$,
$$
|\oph^{\alpha,\beta^\prime}(t^\prime)-\oph^\ab(t)|\leq
|\oph^{\alpha,\beta^\prime}(t^\prime)-\oph^{\alpha,\beta^\prime}(t)|+
|\oph^{\alpha,\beta^\prime}(t)-\oph^\ab(t)|\leq
\frac{|t-t^\prime|}{\beta^\prime}+\frac{|\beta-\beta^\prime|}{\beta^\prime}\,.
$$
(The maximum of $|\oph^{\alpha,\beta^\prime}(t)-\oph^\ab(t)|$ is
taken at $\alpha+\beta$). By induction
$$
|\oph_n^{\alpha,\beta^\prime}(x_0,\ldots,x_{n-1})-
\oph_n^{\alpha,\beta}(x_0,\ldots,x_{n-1})|\leq
|\beta-\beta^\prime|\sum_{j=1}^n(\beta^\prime)^{-j}\,.
$$
Since $\beta^\prime>1$ the sum is convergent. This proves the
first statement. The second statement is proved similarly using
$$
|\oph^\abp(t^\prime)-\oph^\ab(t)|\leq
|\oph^\abp(t^\prime)-\oph^\abp(t)|+|\oph^\abp(t)-\oph^\ab(t)|\leq
\frac{|t-t^\prime|}{\beta^\prime}+\frac{|\alpha-\alpha^\prime|}{\beta^\prime}
$$
which is valid for $\gamma=\gamma^\prime\in [0,1]$ and
$0\leq\alpha^\prime<\alpha\leq 1$. We prove the last statement.
Given $\varepsilon>0$ there exists $n^*$
$$
\oph_{n^*}^{\alpha,1}(\ud{x})\geq
\ophi^{\alpha,1}(\ud{x})-\varepsilon\,.
$$
Since $\beta\mapsto \oph_{n^*}^{\alpha,\beta}(\ud{x})$ is
continuous, there exists $\beta^\prime$ so that for
$1\leq\beta\leq\beta^\p$,
$$
\oph_{n}^{\alpha,\beta}(\ud{x})\geq
\oph_{n^*}^{\alpha,\beta^\p}(\ud{x})\geq
\oph_{n^*}^{\alpha,1}(\ud{x})-\varepsilon\quad\forall n\geq n^*\,.
$$
Hence
$$
\ophi^{\alpha,1}(\ud{x})-2\varepsilon\leq \ophi^\ab(\ud{x})\leq
\ophi^{\alpha,1}(\ud{x})\,.
$$

 \qed

\begin{cor}\label{cor3.2}
Given $\ud{x}$ and $0\leq\alpha^*\leq 1$, let
$$
g_{\alpha^*}(\gamma):=\ophi^{\alpha^*,\beta(\gamma)}(\ud{x})\quad\text{with}\quad
\beta(\gamma):=\gamma-\alpha^*+k-1\,.
$$
For $k\geq 3$ the map $g_{\alpha^*}$ is continuous and
non-increasing on $[0,1]$. If $0<g_{\alpha^*}(\gamma_0)<1$, then
the map is strictly decreasing in a neighborhood of $\gamma_0$. If
$k=2$ then the same statements hold on $[\alpha^*,1]$.
\end{cor}

\begin{cor}\label{cor3.3}
Given $\ud{x}$ and $0<\gamma^*\leq 1$, let
$$
h_{\gamma^*}(\alpha):=\ophi^{\alpha,\beta(\alpha)}(\ud{x})\quad\text{with}\quad
\beta(\alpha):=\gamma^*-\alpha+k-1\,.
$$
For $k\geq 3$ the map $h_{\gamma^*}$ is continuous and
non-increasing on $[0,1]$. If $0<h_{\gamma^*}(\alpha_0)<1$, then
the map is strictly decreasing in a neighborhood of $\alpha_0$. If
$k=2$ then the same statements hold on $[0,\gamma^*)$.
\end{cor}

\begin{pro}\label{pro3.3}
Let $k\geq 2$, $\ud{u}, \ud{v}\in\tA^{\Z_+}$ verifying $u_0=0$ and
$v_0=k-1$ and
$$
\sigma\ud{u}\preceq\ud{v}\quad\text{and}\quad
\ud{u}\preceq\sigma\ud{v}\,.
$$
If $k=2$ we also assume that $\sigma\ud{u}\preceq\sigma\ud{v}$.
Then there exist $\bar{\alpha}\in[0,1]$ and
$\bar{\beta}\in[1,\infty)$ so that $\bar{\gamma}\in[0,1]$. If
$\bar{\beta}>1$, then
$$
\ophi^{\bar{\alpha},\bar{\beta}}(\sigma\ud{u})=\bar{\alpha}\quad\text{and}\quad
\ophi^{\bar{\alpha},\bar{\beta}}(\sigma\ud{v})=\bar{\gamma}\,.
$$
\end{pro}

\prf We consider separately the cases $\sigma\ud{v}=\ud{0}$ and
$\sigma\ud{u}=(k-1)^\infty$ (i.e. $u_j=k-1$ for all $j\geq 1$). If
$\sigma\ud{v}=\ud{0}$, then $\ud{u}=\ud{0}$ and
$\ud{v}=(k-1)\ud{0}$; we set $\bar{\alpha}:=0$ and
$\bar{\beta}:=k-1$ ($\bar{\gamma}=0$). If
$\sigma\ud{u}=(k-1)^\infty$, then $\ud{v}=(k-1)^\infty$ and
$\ud{u}=0(k-1)^\infty$; we set $\bar{\alpha}:=1$ and
$\bar{\beta}:=k$.

From now on we assume that $\ud{0}\prec\sigma\ud{v}$ and
$\sigma\ud{u}\prec (k-1)^\infty$. Set $\alpha_0:=0$ and
$\beta_0:=k$. We consider in details the case $k=2$, so that we
also assume that $\sigma\ud{u}\preceq\sigma\ud{v}$.

\noindent
{\sf Step $1$.\,} Set $\alpha_1:=\alpha_0$ and solve the equation
$$
\ophi^{\alpha_1,\beta}(\sigma\ud{v})=\beta+\alpha_1-k+1\,.
$$
There exists a unique solution, $\beta_1$, such that
$k-1<\beta_1\leq k$. Indeed, the map
$$
G_{\alpha_1}(\gamma):=g_{\alpha_1}(\gamma)-\gamma\quad\text{with}\quad
g_{\alpha_1}(\gamma):=\ophi^{\alpha_1,\beta(\gamma)}(\sigma\ud{v})
\;\,\text{and}\;\,\beta(\gamma):=\gamma-\alpha_1+k-1
$$
is continuous and strictly decreasing on $[\alpha_1,1]$ (see
Corollary \ref{cor3.2}). If $\sigma\ud{v}=(k-1)^\infty$, then
$G_{\alpha_1}(1)=0$ and we set $\beta_1:=k$ and we have
$\gamma_1=1$. If $\sigma\ud{v}\not=(k-1)^\infty$, then there
exists a smallest $j\geq 1$ so that $v_j\leq (k-2)$. Therefore
$\ophi^{\alpha_1,k}(\sigma^j\ud{v})<1$ and
$$
\ophi^{\alpha_1,k}(\sigma\ud{v})=\oph^{\alpha_1,k}_{j-1}\big(v_1,\ldots,v_{j-1}+
\ophi^{\alpha_1,k}(\sigma^j\ud{v})\big)<1\,,
$$
so that $G_{\alpha_1}(1)<0$. On the other hand, since
$\sigma\ud{v}\not=\ud{0}$, $\ophi^{\alpha_1,k-1}(\sigma\ud{v})>0$,
so that $G_{\alpha_1}(0)>0$. There exists a unique
$\gamma_1\in(0,1)$ with $G_{\alpha_1}(\gamma_1)=0$. Define
$\beta_1:=\beta(\gamma_1)=\gamma_1-\alpha_1+k-1$.

\noindent
{\sf Step $2$.\,} Solve in $[0,\gamma_1)$ the equation
$$
\ophi^{\alpha,\beta(\alpha)}(\sigma\ud{u})=\alpha\quad\text{with}\quad\beta(\alpha):=
\gamma_1-\alpha+k-1=\beta_1+\alpha_1-\alpha\,.
$$
If $\sigma\ud{u}=0$, then set $\bar{\alpha}:=0$ and
$\bar{\beta}:=\beta_1$. Let $\sigma\ud{u}\not=0$. There exists a
smallest $j\geq 1$ such that $u_j\geq 1$. This implies that
$\ophi^{\alpha_1,\beta_1}(\sigma^j\ud{u})>0$ and consequently
$$
\ophi^{\alpha_1,\beta(\alpha_1)}(\sigma\ud{u})=
\oph^{\alpha_1,\beta_1}_{j-1}\big(u_1,\ldots,u_{j-1}+
\ophi^{\alpha_1,\beta_1}(\sigma^j\ud{u})\big)>0\,.
$$
Since $\sigma\ud{u}\preceq\sigma\ud{v}$,
$$
0<\ophi^{\alpha_1,\beta_1}(\sigma\ud{u})
\leq\ophi^{\alpha_1,\beta_1}(\sigma\ud{v})=\gamma_1\,.
$$
We have $\gamma_1=1$ only in the case $\sigma\ud{v}=(k-1)^\infty$;
in that case we also have
$\ophi^{\alpha_1,\beta_1}(\sigma\ud{u})<1$. By Corollary
\ref{cor3.3}, for any $\alpha>\alpha_1$ we have
$\ophi^{\alpha_1,\beta_1}(\sigma\ud{u})>\ophi^{\alpha,\beta(\alpha)}(\sigma\ud{u})$.
Therefore, the map
$$
H_{\gamma_1}(\alpha):=h_{\gamma_1}(\alpha)-\alpha\quad\text{with}\quad
h_{\gamma_1}(\alpha):=\ophi^{\alpha,\beta(\alpha)}(\sigma\ud{u})
$$
is continuous and strictly decreasing on $[0,\gamma_1)$,
$H_{\gamma_1}(\alpha_1)>0$ and
$\lim_{\alpha\uparrow\gamma_1}H_{\gamma_1}(\alpha)<0$. There
exists a unique  $\alpha_2\in (\alpha_1,\gamma_1)$ such that
$H_{\gamma_1}(\alpha_2)=0$. Set
$\beta_2:=\gamma_1-\alpha_2+k-1=\alpha_1+\beta_1-\alpha_2$ and
$\gamma_2:=\alpha_2+\beta_2-k+1=\gamma_1$. Since $\alpha_2\in
[0,\gamma_1)$, we have $\beta_2>1$. Hence
\begin{equation}\label{step2}
\alpha_1<\alpha_2<\gamma_1\quad\text{and}\quad
1<\beta_2<\beta_1\quad\text{and}\quad \gamma_2=\gamma_1\,.
\end{equation}
If $\sigma\ud{v}=(k-1)^\infty$, $\gamma_2=1$ and  we set
$\bar{\alpha}:=\alpha_2$ and $\bar{\beta}:=\beta_2$.

\noindent
{\sf Step $3$.\,} From now on $\sigma\ud{u}\not=\ud{0}$ and
$\sigma\ud{v}\not=(k-1)^\infty$. Set $\alpha_3:=\alpha_2$ and
solve in $[\alpha_3,1]$ the equation
$$
\ophi^{\alpha_3,\beta(\gamma)}(\sigma\ud{v})=\gamma\quad\text{with}\quad
\beta(\gamma):=\gamma-\alpha_3+k-1\,.
$$
By Lemma \ref{lem3.3.1}  ($k=2$),
$$
\ophi^{\alpha_3,\beta(\alpha_3)}(\sigma\ud{v})=
\ophi^{\alpha_2,1}(\sigma\ud{v})\geq
\ophi^{\alpha_2,1}(\sigma\ud{u})>\ophi^{\alpha_2,\beta_2}(\sigma\ud{u})=\alpha_2\,,
$$
since $0<\alpha_2<1$. On the other hand by Corollary \ref{cor3.3},
\begin{equation}\label{esti3}
\ophi^{\alpha_3,\beta(\gamma_1)}(\sigma\ud{v})=
\ophi^{\alpha_3,1+\gamma_1-\alpha_3}(\sigma\ud{v})<
\ophi^{\alpha_1,1+\gamma_1-\alpha_1}(\sigma\ud{v})=\ophi^{\alpha_1,\beta_1}(\sigma\ud{v})
=\gamma_1\,,
\end{equation}
since $0<\gamma_1<1$. Therefore, the map $G_{\alpha_3}$ is
continuous and strictly decreasing on $[\alpha_3,1]$,
$G_{\alpha_3}(\alpha_3)>0$ and $G_{\alpha_3}(\gamma_1)<0$. There
exists a unique $\gamma_3\in (\alpha_3,\gamma_1)$ such that
$G_{\alpha_3}(\gamma_3)=0$. Set $\beta_3:=\gamma_3-\alpha_3+k-1$,
so that $\beta_3<\gamma_1-\alpha_2+k-1=\beta_2$. Hence
\begin{equation}\label{step3}
\alpha_3=\alpha_2\quad\text{and}\quad
1<\beta_3<\beta_2\quad\text{and}\quad 0<\gamma_3<\gamma_2<1\,.
\end{equation}

\noindent
{\sf Step $4$.\,} Solve in $[0,\gamma_3)$ the equation
$$
\ophi^{\alpha,\beta(\alpha)}(\sigma\ud{u})=\alpha\quad\text{with}\quad\beta(\alpha):=
\gamma_3-\alpha+k-1=\beta_3+\alpha_3-\alpha\,.
$$
By Lemma \ref{lem3.3.1}
\begin{equation}\label{esti4}
\ophi^{\alpha_3,\beta(\alpha_3)}(\sigma\ud{u})=\ophi^{\alpha_3,\beta_3}(\sigma\ud{u})
>\ophi^{\alpha_3,\beta_2}(\sigma\ud{u})=\ophi^{\alpha_2,\beta_2}(\sigma\ud{u})=\alpha_2\,,
\end{equation}
since $0<\alpha_2<1$. On the other hand,
$$
0<\ophi^{\alpha_3,\beta(\alpha_3)}(\sigma\ud{u})=
\ophi^{\alpha_3,\beta_3}(\sigma\ud{u})\leq\ophi^{\alpha_3,\beta_3}(\sigma\ud{v})
=\gamma_3<1\,.
$$
By Corollary \ref{cor3.3}
$$
\ophi^{\alpha,\beta(\alpha)}(\sigma\ud{u})<\ophi^{\alpha_3,\beta(\alpha_3)}(\sigma\ud{u})
\quad\forall \alpha\in (\alpha_3,\gamma_3)\,.
$$
Therefore, the map
$$
H_{\gamma_3}(\alpha):=h_{\gamma_3}(\alpha)-\alpha\quad\text{with}\quad
h_{\gamma_3}(\alpha):=\ophi^{\alpha,\beta(\alpha)}(\sigma\ud{u})
$$
is continuous and strictly decreasing on $[\alpha_3,\gamma_3)$,
$H_{\gamma_3}(\alpha_3)>0$ and
$\lim_{\alpha\uparrow\gamma_3}H_{\gamma_3}(\alpha) <0$. There
exists a unique  $\alpha_4\in (\alpha_3,\gamma_3)$. Set
$\beta_4:=\gamma_3-\alpha_4+k-1=\alpha_3+\beta_3-\alpha_4$ and
$\gamma_4:=\alpha_4+\beta_4-k+1=\gamma_3$. Hence
\begin{equation}\label{step4}
\alpha_3<\alpha_4<\gamma_3\quad\text{and}\quad
1<\beta_4<\beta_3\quad\text{and}\quad \gamma_4=\gamma_3\,.
\end{equation}
Repeating steps 3 and 4 we get two monotone sequences
$\{\alpha_n\}$ and $\{\beta_n\}$. We set
$\bar{\alpha}:=\lim_{n\ra\infty}\alpha_n$ and
$\bar{\beta}:=\lim_{n\ra\infty}\beta_n$.

We consider briefly the changes which occur when $k\geq 3$. Step
$1$ remains the same. In step $2$ we solve the equation
$H_{\gamma_1}(\alpha)=0$ on $[0,1)$ instead of $[0,\gamma_1)$.
The proof that $H_{\gamma_1}(\alpha_1)>0$ remains the same. We
prove that $\lim_{\alpha\uparrow 1}H_{\gamma_1}(\alpha)<0$.
Corollary \ref{cor3.3} implies that
$$
\gamma_1=\ophi^{\alpha_1,\beta_1}(\sigma\ud{v})=
\ophi^{\alpha_1,\beta(\alpha_1)}(\sigma\ud{v})>
\ophi^{\alpha,\beta(\alpha)}(\sigma\ud{v})\quad
\forall\alpha>\alpha_1\,.
$$
Since $\sigma\ud{u}\preceq \ud{v}$ and $\beta(\alpha_1)=\beta_1$,
$$
\ophi^{\alpha,\beta(\alpha)}(\sigma\ud{u})\leq
\oph^{\alpha,\beta(\alpha)}\big(v_0+\ophi^{\alpha,\beta(\alpha)}(\sigma\ud{v})\big)
\leq \oph^{\alpha_1,\beta(\alpha_1)}
\big(v_0+\ophi^{\alpha,\beta(\alpha)}(\sigma\ud{v})\big)<1\,.
$$
Instead of \eqref{step2} we have
$$
\alpha_1<\alpha_2<1\quad\text{and}\quad
1<\beta_2<\beta_1\quad\text{and}\quad \gamma_2=\gamma_1\,.
$$
Estimate \eqref{esti3} is still valid in step 3 with $k\geq 3$.
Hence $G_{\alpha_3}(\gamma_1)<0$. We solve the equation
$G_{\alpha_3}(\gamma)=0$ on $[0,\gamma_1]$. We have
$$
\ophi^{\alpha_3,\beta(\gamma_1)}(\sigma\ud{u})=
\ophi^{\alpha_2,\beta_2}(\sigma\ud{u})=\alpha_2\,.
$$
By Corollary \ref{cor3.2} we get
$$
\ophi^{\alpha_3,\beta(\gamma)}(\sigma\ud{u})>
\ophi^{\alpha_2,\beta_2}(\sigma\ud{u})=\alpha_2 \quad
\forall\gamma<\gamma_1\,.
$$
Since $\ud{u}\preceq\sigma\ud{v}$,
$$
\ophi^{\alpha_3,\beta(0)}(\sigma\ud{v})\geq
\oph^{\alpha_3,\beta(0)}\big(u_0+\ophi^{\alpha_3,\beta(0)}(\sigma\ud{u})\big)
\geq
\oph^{\alpha_2,\beta(\gamma_1)}\big(u_0+\ophi^{\alpha_3,\beta(0)}(\sigma\ud{u})\big)>0\,.
$$
Estimate \eqref{esti4} is still valid in step 4 so that
$H_{\gamma_3}(\alpha_3)>0$. Corollary \ref{cor3.3} implies that
$$
\gamma_3=\ophi^{\alpha_3,\beta_3}(\sigma\ud{v})=
\ophi^{\alpha_3,\beta(\alpha_3)}(\sigma\ud{v})>
\ophi^{\alpha,\beta(\alpha)}(\sigma\ud{v})\quad
\forall\alpha>\alpha_3\,.
$$
Therefore
$$
\ophi^{\alpha_3,\beta(\alpha_3)}(\sigma\ud{u})\leq
\oph^{\alpha,\beta(\alpha)}\big(v_0+\ophi^{\alpha,\beta(\alpha)}(\sigma\ud{v})\big)
\leq \oph^{\alpha_3,\beta(\alpha_3)}
\big(v_0+\ophi^{\alpha,\beta(\alpha)}(\sigma\ud{v})\big)<1\,.
$$
Instead of \eqref{step4} we have
$$
\alpha_3<\alpha_4<1\quad\text{and}\quad
1<\beta_4<\beta_3\quad\text{and}\quad \gamma_4=\gamma_3\,.
$$

Assume that $\bar{\beta}>1$. Then $1<\bar{\beta}\leq\beta_n$ for
all $n$. We have
$$
\ophi^{\alpha_n,\beta_n}(\sigma\ud{v})=\gamma_n\,,\quad\text{$n$
odd}
$$
and
$$
\ophi^{\alpha_n,\beta_n}(\sigma\ud{u})=\alpha_n\,,\quad\text{$n$
even}\,.
$$
Let $\bar{\gamma}=\bar{\alpha}+\bar{\beta}-k+1$. For $n$ odd, let
$\beta_n^*:=\bar{\gamma}-\alpha_n+k-1$; using Lemma \ref{lem3.3.2}
we get
\begin{align*}
|\ophi^{\bar{\alpha},\bar{\beta}}(\sigma\ud{v})-\bar{\gamma}|&\leq
|\ophi^{\bar{\alpha},\bar{\beta}}(\sigma\ud{v})-
\ophi^{\alpha_n,\beta_n^*}(\sigma\ud{v})|+
|\ophi^{\alpha_n,\beta_n^*}(\sigma\ud{v})-
\ophi^{\alpha_n,\beta_n}(\sigma\ud{v})|+|\gamma_n-\gamma|\\
&\leq\frac{1}{\bar{\beta}-1}(2|\bar{\alpha}-\alpha_n|+
|\bar{\beta}-\beta_n|)+|\gamma_n-\gamma|\,,
\end{align*}
since  $\beta_n^*=\bar{\beta}+\bar{\alpha}-\alpha_n$. Letting $n$
going to infinity we get
$\ophi^{\bar{\alpha},\bar{\beta}}(\sigma\ud{v})=\bar{\gamma}$.
Similarly we prove
$\ophi^{\bar{\alpha},\bar{\beta}}(\sigma\ud{v})=\bar{\alpha}$.
\qed

\begin{cor}\label{cor3.4}
Suppose that $(\ud{u},\ud{v})$, respectively
$(\ud{u}^\p,\ud{v}^\p)$, verify the hypothesis of Proposition
\ref{pro3.3} with $k\geq 2$, respectively with $k^\prime\geq 2$.
If $k\geq k^\prime$, $\ud{u}\preceq\ud{u}^\p$ and
$\ud{v}^\p\preceq\ud{v}$, then $\bar{\beta}^\p\leq\bar{\beta}$
and $\bar{\alpha}^\p\geq \bar{\alpha}$.
\end{cor}

\prf We consider the case $k=k^\p$, whence
$\sigma\ud{v}^\p\preceq\sigma\ud{v}$. From the proof of
Proposition \ref{pro3.3} we get $\gamma^\p_1\leq\gamma_1$ and
$\alpha_1^\p\geq\alpha_1$. Suppose that $\gamma^\p_j\leq\gamma_j$
and $\alpha_j^\p\geq\alpha_j$ for $j=1,\ldots,n$. If $n$ is even,
then $\alpha_{n+1}^\p=\alpha_{n}^\p$ and
$\alpha_{n+1}=\alpha_{n}$. We prove that
$\gamma_{n+1}^\p\leq\gamma_{n+1}$. We have
$$
\gamma_{n+1}^\p=\ophi^{\alpha_{n+1}^\p,\beta(\gamma_{n+1}^\p)}(\sigma\ud{v}^\p)
\leq
\ophi^{\alpha_{n+1}^\p,\beta(\gamma_{n+1}^\p)}(\sigma\ud{v})\leq
\ophi^{\alpha_{n+1},\beta(\gamma_{n+1}^\p)}(\sigma\ud{v})\implies
\gamma_{n+1}\geq\gamma_{n+1}^\p\,.
$$
If $n$ is odd, then  $\gamma_{n+1}^\p=\gamma_{n}^\p$ and
$\gamma_{n+1}=\gamma_{n}$. We prove that
$\alpha_{n+1}^\p\geq\alpha_{n+1}$. We have
\begin{align*}
\alpha_{n+1}&=\ophi^{\alpha_{n+1},\beta(\alpha_{n+1})}(\sigma\ud{u})
\leq \ophi^{\alpha_{n+1},\beta(\alpha_{n+1})}(\sigma\ud{u}^\p)=
\ophi^{\alpha_{n+1},\gamma_{n+1}-\alpha_{n+1}+k-1)}(\sigma\ud{u}^\p)\\
& \leq
\ophi^{\alpha_{n+1},\gamma_{n+1}^\p-\alpha_{n+1}+k-1)}(\sigma\ud{u}^\p)
\implies \alpha_{n+1}^\p\geq\alpha_{n+1}\,.
\end{align*}
\qed

We state a uniqueness  result. The proof uses Theorem
\ref{thm3.1}.

\begin{pro}\label{prounicity}
Let $k\geq 2$, $\ud{u}, \ud{v}\in\tA^{\Z_+}$, $u_0=0$ and
$v_0=k-1$, and assume that \eqref{condition} holds. Then there is
at most one solution $(\ab)\in[0,1]\times[1,\infty)$ for the
equations
$$
\ophi^{\alpha,\beta}(\sigma\ud{u})=\alpha\quad\text{and}\quad
\ophi^{\alpha,\beta}(\sigma\ud{v})=\gamma\,.
$$
\end{pro}

\prf Assume that there are two solutions $(\alpha_1,\beta_1)$ and
$(\alpha_2,\beta_2)$ with $\beta_1\leq \beta_2$. If
$\alpha_2>\alpha_1$, then
$$
\alpha_2-\alpha_1=\ophi^{\alpha_2,\beta_2}(\sigma\ud{u})-
\ophi^{\alpha_1,\beta_1}(\sigma\ud{u})\leq 0\,,
$$
which is impossible. Therefore $\alpha_2\leq\alpha_1$. If
$\beta_1=\beta_2$, then
$$
0\geq\alpha_2-\alpha_1=\ophi^{\alpha_2,\beta_2}(\sigma\ud{u})-
\ophi^{\alpha_1,\beta_1}(\sigma\ud{u})\geq 0\,,
$$
which implies $\alpha_2=\alpha_1$. Therefore we assume that
$\alpha_2\leq\alpha_1$ and $\beta_1<\beta_2$. However, Theorem
\ref{thm3.1} implies that
$$
\log_2\beta_1=h(\BSigma(\ud{u},\ud{v}))=\log_2\beta_2\,,
$$
which is impossible. \qed

\subsection{Computation of the topological entropy of
$\BSigma(\ud{u},\ud{v})$}\label{topological}

We compute the entropy of the shift space $\BSigma(\ud{u},\ud{v})$
where $\ud{u}$ and $\ud{v}$  is a pair of strings verifying
$u_0=0$, $v_0=k-1$ and \eqref{condition}. The main result is
Theorem \ref{thm3.1}. The idea for computing the topological
entropy is to compute  $\bar{\alpha}$ and $\bar{\beta}$ by the
algorithm of section \ref{subsectionalgo} and to use the fact that
$h\big(\BSigma(\ud{u}^\abb,\ud{v}^\abb)\big)=\log_2\bar{\beta}$
(see e.g. \cite{Ho1}). The most difficult case is when $\ud{u}$
and $\ud{v}$ are both periodic. Assume that the string
$\ud{u}:=\ud{a}^\infty$ has minimal period $p$, $|\ud{a}|=p$, and
that the string $\ud{v}:=\ud{b}^\infty$ has minimal period $q$,
$|\ud{b}|=q$. If $a_0=a_{p-1}=0$, then $\ud{u}=\ud{0}$ and $p=1$.
Indeed, if $a_0=a_{p-1}=0$, then
$\ud{a}\ud{a}=(\tp\ud{a})00(\s\ud{a})$; the result follows from
\eqref{condition}. Similarly, if $b_0=b_{q-1}=k-1$, then
$\ud{v}=(k-1)^\infty$ and $q=1$. These cases are similar to the
case when  only one of the strings $\ud{u}$ and $\ud{v}$ is
periodic and are simpler than the generic case of two periodic
strings, which we treat in details.

The setting for subsection \ref{topological} is the following one.
The string $\ud{u}:=\ud{a}^\infty$ has minimal period $p\geq 2$
with $u_0=0$, or $\ud{u}=(1)^\infty$. The string
$\ud{v}:=\ud{b}^\infty$ has minimal period $q\geq 2$ with
$v_0=k-1$, or $\ud{v}=(k-2)^\infty$. We also consider the strings
$\ud{u}^*=\ud{a}^\p\ud{b}^\infty$ and
$\ud{v}^*=\ud{b}^\p\ud{a}^\infty$ with
$\ud{a}^\p=\tp\ud{a}(u_{p-1}-1)$ and
$\ud{b}^\p=\tp\ud{b}(v_{q-1}+1)$. We write
$\BSigma\equiv\BSigma(\ud{u},\ud{v})$, $\BSigma^*\equiv
\BSigma(\ud{u}^*,\ud{v}^*)$, $\cG\equiv \cG(\ud{u},\ud{v})$ and
$\cG^*\equiv \cG(\ud{u}^*,\ud{v}^*)$. The main point is to prove
that $h(\cG)=h(\cG^*)$ by comparing the follower-set graphs $\cG$
and $\cG^*$.

\begin{lem}\label{lemgraph}
1) In the above setting the vertices of the graph $\cG$ are
$\cF_\epsilon$, $\cF_{\ud{w}}$ with $\ud{w}$ a prefix of
$\tp\ud{a}$ or of $\tp\ud{b}$, $\tp\ud{a}$ and $\tp\ud{b}$
included. \\
2) Let $r:=|v(\tp\ud{a})|$. If $u_{p-1}\not=v_r$, then
$\cF_{\ud{a}}=\cF_\epsilon$ and there is an edge labeled by
$u_{p-1}$ from $\cF_{\tp\ud{a}}$ to $\cF_\epsilon$. If
$u_{p-1}=v_r$, then $\cF_{\ud{a}}=\cF_{v(\tp\ud{a})v_r}$ and there
is a single edge, labeled by $u_{p-1}=v_r$, from $\cF_{\tp\ud{a}}$
to $\cF_{v(\tp\ud{a})v_r}$. If $k=2$ the first possibility is
excluded.\\
3) Let $s:=|u(\tp\ud{b})|$. If $v_{q-1}\not=u_s$, then
$\cF_{\ud{b}}=\cF_\epsilon$ and there is an edge labeled by
$v_{q-1}$ from $\cF_{\tp\ud{b}}$ to $\cF_\epsilon$. If
$v_{q-1}=u_s$, then $\cF_{\ud{b}}=\cF_{u(\tp\ud{b})u_s}$ and there
is a single edge, labeled by $v_{q-1}=u_s$, from $\cF_{\tp\ud{b}}$
to $\cF_{u(\tp\ud{b})u_s}$. If $k=2$ the first possibility is
excluded.
\end{lem}

\prf Suppose that $\ud{w}$ and $\ud{w}\ud{w}^\p$ are two prefixes
of $\tp\ud{a}$. We show that
$\cF_{\ud{w}}\not=\cF_{\ud{w}\ud{w}^\p}$. Write
$\ud{u}=\ud{w}\ud{x}$ and $\ud{u}=\ud{w}\ud{w}^\p\ud{y}$ and
suppose that $\cF_{\ud{w}}=\cF_{\ud{w}\ud{w}^\p}$. Then (see
\eqref{3.2.3}) $\ud{x}=\ud{y}=\sigma^p\ud{u}$, so that
$\ud{u}=(\ud{w}^\prime)^\infty$, contradicting the minimality of
the period $p$. Consider the vertex $\cF_{\tp\ud{a}}$ of $\cG$.
We have
$$
\cF_{\tp\ud{a}}=\{\ud{x}\in\BSigma\colon
\sigma^{p-1}\ud{u}\preceq\ud{x}\preceq
\sigma^{r}\ud{v}\}\quad\text{where $r=|v(\tp\ud{a})|$}\,.
$$
Let $\ud{d}$ be the prefix of $\ud{v}$ of length $r+1$, so that
$\tp\ud{d}=v(\tp\ud{a})$.  If $u_{p-1}\not=v_{r}$, then there are
an edge labeled by $u_{p-1}$ from $\cF_{\tp\ud{a}}$ to
$\cF_{\ud{a}}=\cF_{\epsilon}$ (since $\sigma^p\ud{u}=\ud{u}$) and
an edge labeled by $v_{r}$ from $\cF_{\tp\ud{a}}$ to
$\cF_{\ud{d}}$. There may be other labeled edges from
$\cF_{\tp\ud{a}}$ to $\cF_\epsilon$ (see Lemma \ref{lem3.2.2}). If
$u_{p-1}=v_{r}$, then there is a single out-going  edge labeled by
$u_{p-1}$ from $\cF_{\tp\ud{a}}$ to $\cF_{\ud{a}}$ and
$v(\ud{a})=\ud{d}$. We prove that $\cF_{\ud{a}}=\cF_{\ud{d}}$. If
$u(\ud{d})=\epsilon$, the result is true, since in that case
$$
\cF_{\ud{d}}=\{\ud{x}\in\BSigma\colon \ud{u}\preceq\ud{x}\preceq
\sigma^{r+1}\ud{v}\}=\{\ud{x}\in\BSigma\colon
\sigma^p\ud{u}\preceq\ud{x}\preceq
\sigma^{r+1}\ud{v}\}=\cF_{\ud{a}}\,.
$$
We exclude the possibility $u(\ud{d})\not=\epsilon$.  Suppose that
$\ud{w}:=u(\ud{d})$ is non-trivial ($|u(\ud{d})|<p$). We can write
$\ud{a}=\ud{a}^{\p\p}\ud{w}$ and $\ud{a}=\ud{w}\widehat{\ud{a}}$
since $\ud{w}$ is a prefix of $\ud{u}$, and consequently
$\ud{a}\ud{a}=\ud{a}^{\p\p}\ud{w}\ud{w}\widehat{\ud{a}}$. From
Lemma \ref{lem3.2.0} we conclude that $\ud{w}\ud{w}$ is a prefix
of $\ud{u}$, so that
$\ud{a}\ud{u}=\ud{a}^{\p\p}\ud{w}\ud{w}\ud{w}\cdots$, proving that
$\ud{u}$ has period $|\ud{w}|$, contradicting the hypothesis that
$p$ is the minimal period of $\ud{u}$. If $k=2$ the first
possibility is excluded because $u_{p-1}\not=0$ and we have
$u_{p-1}\preceq v_r$ by $\sigma^{p-r}\ud{u}\preceq\ud{v}$. The
discussion concerning the vertex $\cF_{\tp\ud{b}}$ is similar.
\qed

\begin{pro}\label{pro3.6}
Consider the above setting. If $h(\BSigma)>0$, then
$h(\BSigma)=h(\BSigma^*)$.
\end{pro}

\prf Consider the vertex $\cF_{\tp\ud{a}}$ of $\cG^*$. In that
case $(u_{p-1}-1)\not=v_{r}$ so that we have an additional edge
labeled by $u_{p-1}-1$ from $\cF_{\tp\ud{a}}$ to
$\cF_{\ud{a}^\p}$ (see proof of Lemma \ref{lemgraph}), otherwise
all out-going edges from $\cF_{\tp\ud{a}}$, which are present in
the graph $\cG$, are also present in $\cG^*$. Let $v^*(\ud{w})$
be the longest suffix of $\ud{w}$, which is a prefix of
$\ud{v}^*$. Then
$$
\cF_{\ud{a}^\p}=\{\ud{x}\in\BSigma^*\colon
\sigma^{p}\ud{u}^*\preceq\ud{x}\preceq \ud{v}^*\}=
\{\ud{x}\in\BSigma^*\colon \ud{v}\preceq\ud{x}\preceq
\ud{v}^*\}\,.
$$
Similarly, there is an additional edge labeled by $v_{q-1}+1$
from $\cF_{\tp\ud{b}}$ to $\cF_{\ud{b}^\p}$. Let $u^*(\ud{w})$ be
the longest suffix of $\ud{w}$, which is a prefix of $\ud{u}^*$.
Then
$$
\cF_{\ud{b}^\p}=\{\ud{x}\in\BSigma^*\colon
\ud{u}^*\preceq\ud{x}\preceq \sigma^q\ud{v}^*\}=
\{\ud{x}\in\BSigma^*\colon \ud{u}^*\preceq\ud{x}\preceq
\ud{u}\}\,.
$$
The structure of the graph $\cG^*$ is very simple from the
vertices $\cF_{\ud{a}^\p}$ and $\cF_{\ud{b}^\p}$. There is a
single  out-going edge from $\cF_{\ud{a}^\p}$ to
$\cF_{\ud{a}^{\p}v_0}$, from $\cF_{\ud{a}^{\p}v_0}$ to
$\cF_{\ud{a}^{\p}v_0v_1}$ and so on, until we reach the vertex
$\cF_{\ud{a}^{\p}\tp\ud{b}}$. From that vertex there are an
out-going edge labeled by $v_{q-1}$ to $\cF_{\ud{a}^\p}$ and an
out-going edge labeled by $v_{q-1}+1$ to $\cF_{\ud{b}^\p}$.
Similarly, there is a single out-going edge from $\cF_{\ud{b}^\p}$
to $\cF_{\ud{b}^{\p}u_0}$, from $\cF_{\ud{b}^{\p}u_0}$ to
$\cF_{\ud{b}^{\p}u_0u_1}$ and so on, until we reach the vertex
$\cF_{\ud{b}^{\p}\tp\ud{a}}$. From that vertex there are an
out-going  edge labeled by $u_{p-1}$ to $\cF_{\ud{b}^\p}$ and an
out-going  edge labeled by $u_{p-1}-1$ to $\cF_{\ud{a}^\p}$. Let
us denote that part of $\cG^*$ by $\cG^*\backslash\cG$. This
subgraph is strongly connected. The graph $\cG^*$ consists of the
union of $\cG$ and $\cG^*\backslash\cG$ with the addition of the
two edges from $\cF_{\tp\ud{a}}$ to $\cF_{\ud{a}^\p}$ and
$\cF_{\tp\ud{b}}$ to $\cF_{\ud{b}^\p}$. Using Theorem 1.7 of
\cite{BGMY} it easy to compute the entropy of the subgraph
$\cG^*\backslash\cG$ (use as rome $\{\cF_{\ud{a}^\p},
\cF_{\ud{b}^\p}\}$). It is the largest root, say $\lambda^*$, of
the equation
$$
\lambda^{-q}+\lambda^{-p}-1=0\,.
$$
Hence $\lambda^*$ is equal to the entropy of a graph with two
cycles of periods $p$ and $q$, rooted at a  common point. To
prove Proposition \ref{pro3.6} it sufficient to exhibit a
subgraph of $\cG$ which has an entropy larger or equal to that of
$\cG^*\backslash\cG$.

If $k\geq 4$, then there is a subgraph with two cycles of length
$1$ rooted at $\cF_\epsilon$. Hence $h(\cG)\geq
\log_22>\lambda^*$.
 If $\cF_{\ud{a}}=\cF_\epsilon$ or
$\cF_{\ud{b}}=\cF_\epsilon$, which could happen only for $k\geq
3$ (see Lemma \ref{lemgraph}), then there is a subgraph of $\cG$
consisting of two cycles rooted at $\cF_\epsilon$, one of length
$p$ or  of length $q$ and another one of length $1$. This also
implies that $h(\cG)\geq\lambda^*$. Since the minimal periods of
$\ud{u}$ and $\ud{v}$ are $p$ and $q$, it is impossible that
$\cF_{\ud{w}}=\cF_\epsilon$ for $\ud{w}$ a non trivial prefix of
$\tp\ud{a}$ or $\tp\ud{b}$. Therefore we assume that $k\leq 3$,
$\cF_{\ud{a}}\not=\cF_\epsilon$ and
$\cF_{\ud{b}}\not=\cF_{\epsilon}$.

Let $\cH$ be a  strongly connected component of $\cG$ which has
strictly positive entropy. If $\cF_\epsilon$ is a vertex of $\cH$,
which happens only if $k=3$, then we conclude as above that
$h(\cG)\geq\lambda^*$. Hence, we assume that $\cF_\epsilon$ is not
a vertex of $\cH$. The vertices of $\cH$ are indexed by prefixes
of $\tp\ud{a}$ and $\tp\ud{b}$. Let $\cF_{\ud{c}}$  be the vertex
of $\cH$ with $\ud{c}$ a prefix of $\ud{u}$ and $|\ud{c}|$
minimal; similarly, let $\cF_{\ud{d}}$ be the vertex of $\cH$ with
$\ud{d}$ a prefix of $\ud{v}$ and $|\ud{d}|$ minimal. By our
assumptions $r:=|\ud{c}|\geq 1$ and $s:=|\ud{d}|\geq 1$. The
following argument is a simplified adaptation of the proof of
Lemma 3 in \cite{Ho3}. The core of the argument is the content of
the Scholium \ref{sch}. Consider the $v$-parsing of $\ud{a}$ from
the prefix $\tp\ud{c}$, and the $u$-parsing of $\ud{b}$ from the
prefix $\tp\ud{d}$,
$$
\ud{a}=(\tp\ud{c})\ud{a}^1\cdots\ud{a}^k\quad\text{and}\quad
\ud{b}=(\tp\ud{d})\ud{b}^1\cdots\ud{b}^\ell\,.
$$
(From $\tp\ud{c}$ the $v$-parsing of $\ud{a}$ does not depend on
$\tp\ud{c}$ since there is an in-going edge at $\cF_{\ud{c}}$.)

 We claim that there are an edge from
$\cF_{(\tp\ud{c})\ud{a}^1}$ to $\cF_{\ud{d}}$ and an edge from
$\cF_{(\tp\ud{d})\ud{b}^1}$ to $\cF_{\ud{c}}$. Suppose that this
is not the case, for example, there is an edge from
$\cF_{(\tp\ud{c})\ud{a}^1\cdots \ud{a}^j}$ to $\cF_{\ud{d}}$, but
no edge from $\cF_{(\tp\ud{c})\ud{a}^1\cdots \ud{a}^i}$ to
$\cF_{\ud{d}}$, $1\leq i<j$. This implies that
$v(\ud{a}^j)=\s\ud{a}^j=\tp(\ud{d})$ and
$(\tp\ud{c})\ud{a}^1\cdots \ud{a}^jf^\prime$ is a prefix of
$\ud{u}$ with $f^\p\prec f$ and $f$ defined by
$\ud{d}=(\tp\ud{d})f$. On the other hand there exists an edge from
$\cF_{(\tp\ud{c})\ud{a}^1\cdots \ud{a}^{j-1}}$ to
$\cF_{v((\tp\ud{c})\ud{a}^1\cdots \ud{a}^{j-1})*}=
\cF_{v(\ud{a}^{j-1})*}$ with $*$ some letter of $\tA$ and
$v(\ud{a}^{j-1})*\not=\ud{d}$ by hypothesis. Let $e$ be the first
letter of $\ud{a}^j$. Then $*=(e+1)$ since we assume that
$\cF_\epsilon$ is not a vertex of $\cH$ and consequently there
are only two out-going edges from $\cF_{(\tp\ud{c})\ud{a}^1\cdots
\ud{a}^{j-1}}$. There exists an edge from $\cF_{v(\ud{a}^{j-1})}$
to $\cF_{u(v(\ud{a}^{j-1}))*}$, where $*$ is some letter of $\tA$
(see Scholium \ref{sch}). Again, since $\cF_\epsilon$ is not a
vertex of $\cH$ we must have $*=e$. Either
$u(v(\ud{a}^{j-1}))e=\ud{c}$ or $u(v(\ud{a}^{j-1}))e\not=\ud{c}$.
In the latter case, by the same reasoning,  there exists an edge
from $\cF_{u(v(\ud{a}^{j-1}))}$ to
$\cF_{v(u(v(\ud{a}^{j-1})))(e+1)}$ and
$v(u(v(\ud{a}^{j-1})))(e+1)\not=\ud{d}$ by hypothesis; there
exists also an edge from $\cF_{v(u(v(\ud{a}^{j-1})))}$ to
$\cF_{u(v(u(v(\ud{a}^{j-1}))))e}$. After a finite number of steps
we get
$$
u(\cdots v(u(v(\ud{a}^{j-1}))))e=\ud{c}\,.
$$
This implies that $\tp\ud{c}$ is a suffix of $\ud{a}^{j-1}$, and
the last letter of $\ud{c}$ (or the first letter of $\ud{a}^1$)
is $e$. Hence $\ud{a}^1=e\ud{d}\cdots$.  If we write
$\ud{a}^{j-1}=\ud{g}(\tp\ud{c})$ we have
$$
(\tp\ud{c})\ud{a}^1=\ud{c}\ud{d}\cdots=\ud{c}(\tp\ud{d})f\cdots
\quad\text{and}\quad\ud{a}^{j-1}\ud{a}^j
f=\ud{g}(\tp\ud{c})e(\tp\ud{d})f^\prime=\ud{g}\ud{c}(\tp\ud{d})f^\prime\,.
$$
We get a contradiction with \eqref{condition} since
$\ud{c}(\tp\ud{d})f^\p\prec\ud{c}(\tp\ud{d})f$.

Consider the smallest strongly connected subgraph $\cH^\p$ of
$\cH$ which contains the vertices $\cF_{\ud{c}}$,
$\cF_{(\tp\ud{c})\ud{a}^1}$, $\cF_{\ud{d}}$ and
$\cF_{(\tp\ud{d})\ud{b}^1}$. Since $\cH$ has strictly positive
entropy, there exists at least one edge from some other vertex $A$
of $\cH$ to $\cF_{\ud{c}}$ or $\cF_{\ud{d}}$, say $\cF_{\ud{c}}$.
Define $\cG^\p$ as the smallest strongly connected subgraph of
$\cH$, which contains $\cH^\p$ and $A$. This graph has two cycles:
one passing through the vertices $\cF_{\ud{c}}$,
$\cF_{(\tp\ud{c})\ud{a}^1}$, $\cF_{\ud{d}}$,
$\cF_{(\tp\ud{d})\ud{b}^1}$ and $\cF_{\ud{c}}$, the other one
passing through the vertices $\cF_{\ud{c}}$,
$\cF_{(\tp\ud{c})\ud{a}^1}$, $\cF_{\ud{d}}$,
$\cF_{(\tp\ud{d})\ud{b}^1}$, $A$ and $\cF_{\ud{c}}$. The first
cycle has length $|\ud{a}^1|+|\ud{b}^1|$, and the second cycle has
length $|\ud{a}^1|+|\ud{b}^1|+\cdots+|\ud{b}^j|$ if
$A=\cF_{\tp(\ud{d})\ud{b}^1\cdots\ud{b}^j}$. We also have
$$
|\ud{c}|=|\ud{b}^1|=|\ud{b}^j|\quad\text{and}\quad
|\ud{a}^1|=|\ud{d}|\,.
$$
Therefore one cycle has period
$$
|\ud{a}^1|+|\ud{b}^1|\leq |\ud{a}^1|+|\ud{c}|\leq p\,,
$$
and the other one has period
$$
|\ud{d}|+|\ud{b}^1|+\cdots+|\ud{b}^j|\leq q\,.
$$
\qed

\begin{thm}\label{thm3.1}
Let $k\geq 2$ and let $\ud{u}\in\tA^{\Z_+}$ and
$\ud{v}\in\tA^{\Z_+}$, such that $u_0=0$, $v_0=k-1$ and
$$
\ud{u}\preceq\sigma^n\ud{u}\preceq \ud{v}\quad \forall\,n\geq
0\quad\text{and}\quad \ud{u}\preceq \sigma^n\ud{v}\preceq
\ud{v}\quad \forall\,n\geq 0\,.
$$
If $k=2$ we also assume that $\sigma\ud{u}\preceq\sigma\ud{v}$.
Let $\bar{\alpha}$ and $\bar{\beta}$ be the two real numbers
defined by the algorithm of Proposition \ref{pro3.3}. Then
$$
h(\BSigma(\ud{u},\ud{v}))=\log_2\bar{\beta}\,.
$$
If $k=2$ and $\sigma\ud{v}\prec \sigma\ud{u}$, then
$h(\BSigma(\ud{u},\ud{v}))=0$.
\end{thm}

\prf Let $\bar{\beta}>1$. By Propositions \ref{pro3.3} and
\ref{pro3.4ter} we have
$$
\BSigma(\ud{u}^\abb,\ud{v}^\abb)\subset
\BSigma(\ud{u},\ud{v})\subset\BSigma(\ud{u}^{\abb}_*,\ud{v}^{\abb}_*)
\,.
$$
From Proposition \ref{pro3.6} we get
$$
h(\BSigma(\ud{u}^\abb,\ud{v}^\abb))=
h(\BSigma(\ud{u}^{\abb}_*,\ud{v}^{\abb}_*))=\log_2\bar{\beta}\,.
$$
Let $\lim_n\alpha_n=\bar{\alpha}$ and
$\lim_n\beta_n=\bar{\beta}=1$.  We have $\alpha_n<1$ and
$\beta_n>1$ (see proof of Proposition \ref{pro3.3}). Let
$$
\ud{u}^n:=\ud{u}^{\alpha_n,\beta_n}_*\quad\text{and}\quad
\ud{v}^n:=\ud{v}^{\alpha_n,\beta_n}_*\,.
$$
By Proposition \ref{pro3.4ter} point 3,
$$
\ud{v}^{\alpha_1,\beta_1}\preceq \ud{v}\preceq \ud{v}^1\,.
$$
By monotonicity,
$$
\ophi^{\alpha_2,\beta_2}(\sigma\ud{v}^1)\leq
\ophi^{\alpha_1,\beta_1}(\sigma\ud{v}^1)=\gamma_1=\gamma_2
\ophi^{\alpha_2,\beta_2}(\sigma\ud{v}^2)\,.
$$
Therefore $\ud{v}^1\preceq\ud{v}^2$ ($v^1_0=v^2_0$) and by
Proposition \ref{pro3.4ter} point 2,
$$
\ud{u}^2\preceq \ud{u}\preceq
\ud{u}^{\alpha_2,\beta_2}\quad\text{and}\quad
\ud{v}\preceq\ud{v}^2\,.
$$
By monotonicity,
$$
\ophi^{\alpha_3,\beta_3}(\sigma\ud{u}^3)=\alpha_3=\alpha_2=
\ophi^{\alpha_2,\beta_2}(\sigma\ud{u}^2)\leq
\ophi^{\alpha_3,\beta_3}(\sigma\ud{u}^2)\,.
$$
Therefore $\ud{u}^3\preceq\ud{u}^2$ and
$$
\ud{u}^3\preceq \ud{u}\quad\text{and}\quad
\ud{v}^{\alpha_3,\beta_3}\preceq\ud{v}\preceq\ud{v}^3\,.
$$
Iterating this argument we conclude that
$$
\ud{u}^n\preceq\ud{u}\quad\text{and}\quad \ud{v}\preceq\ud{v}^n\,.
$$
These inequalities imply
$$
h(\BSigma(\ud{u},\ud{v}))\leq
h(\BSigma(\ud{u}^n,\ud{v}^n))=\log_2\beta_n\ra 0\quad\text{for
$n\ra\infty$}\,.
$$

Finally let $k=2$ and $\sigma\ud{v}\prec\sigma\ud{u}$. If
$\sigma\ud{u}=(1)^\infty$, then $\ud{v}_j=0$ for a single value
of $j$, so that $h(\BSigma(\ud{u},\ud{v}))=0$. Suppose that
$\sigma\ud{u}\not=(1)^\infty$ and fix any $\beta>1$. The function
$\alpha\mapsto\ophi^{\alpha,\beta}(\sigma\ud{u})$ is continuous
and decreasing since $\oph^{\alpha,\beta}$ dominates
$\oph^{\alpha^\p,\beta}$ if $\alpha<\alpha^\p$. There exists
$\alpha\in(0,1)$ such that
$\ophi^{\alpha,\beta}(\sigma\ud{u})=\alpha$. If $\ud{v}_0<
\ud{v}^\ab_0$, then $\ud{v}\prec \ud{v}^\ab$ and
$\BSigma(\ud{u},\ud{v})\subset\BSigma(\ud{u},\ud{v}^\ab)$, whence
$h(\BSigma(\ud{u},\ud{v}))\leq\log_2\beta$. If $\ud{v}_0=
\ud{v}^\ab_0=1$, then
$$
\ophi^\ab(\sigma\ud{v})\leq \ophi^\ab(\sigma\ud{u})=\alpha<\gamma=
\ophi^\ab(\sigma\ud{v}^\ab)\,.
$$
The map $\ophi^\ab$ is continuous and non-decreasing on
$\tA^{\Z_+}$ so that $\sigma\ud{v}\prec \sigma\ud{v}^\ab$, whence
$\ud{v}\prec\ud{v}^\ab$ and
$h(\BSigma(\ud{u},\ud{v}))\leq\log_2\beta$. Since $\beta>1$ is
arbitrary, $h(\BSigma(\ud{u},\ud{v}))=0$.
 \qed

\section{Inverse problem for $\beta x+\alpha \mod 1$}\label{section4}
\setcounter{equation}{0}

In this section we solve the inverse problem for $\beta
x+\alpha\mod 1$, namely the question:  {\it given two strings
$\ud{u}$ and $\ud{v}$ verifying
\begin{equation}\label{4.1}
\ud{u}\preceq\sigma^n\ud{u}\prec\ud{v} \quad\text{and}\quad
\ud{u}\prec\sigma^n\ud{v}\preceq\ud{v} \quad\forall n\geq0\,,
\end{equation}
can we find $\alpha\in [0,1)$ and $\beta\in(1,\infty)$ so that
$\ud{u}=\ud{u}^{\alpha,\beta}$ and
$\ud{v}=\ud{v}^{\alpha,\beta}$? }

\begin{pro}\label{pro3.5}
Let the $\varphi$-expansion be valid. Let $\ud{u}$ be a solution
of \eqref{eqalpha} and  $\ud{v}$ a solution of \eqref{eqbeta}. If
\eqref{4.1} holds, then
$$
\ud{u}^\ab=\ud{u}\iff \forall n\geq
0\colon\;\ophi^\ab(\sigma^n\ud{u})<1 \,\iff\, \forall n\geq
0\colon\;\ophi^\ab(\sigma^n\ud{v})>0\iff \ud{v}^\ab=\ud{v}\,.
$$
\end{pro}

\prf  The $\varphi$-expansion is valid, so that \eqref{validity1}
is true,
$$
\forall n\geq
0\colon\;\ophi^\ab(\sigma^n\ud{u}^\ab)=T^n_\ab(0)<1\,.
$$
Proposition \ref{pro3.4} and Proposition \ref{pro3.4ter} point 2
imply
$$
\ud{u}=\ud{u}^\ab\,\iff\,\forall n\geq
0\colon\;\ophi^\ab(\sigma^n\ud{u})<1\,.
$$
Similarly
$$
\ud{v}=\ud{v}^\ab\,\iff\,\forall n\geq
0\colon\;\ophi^\ab(\sigma^n\ud{v})>0\,.
$$
Let $\ud{x}\prec\ud{x}^\p$, $\ud{x}, \ud{x}^\p\in
\BSigma(\ud{u},\ud{v})$. Let $\ell:=\min\{m\geq 0\colon
x_m\not=x_m^\p\}$. Then
$$
\ophi^\ab(\ud{x})=\ophi^\ab(\ud{x}^\p)\implies
\ophi^\ab(\sigma^{\ell+1}\ud{x})=1\quad\text{and}\quad
\ophi^\ab(\sigma^{\ell+1}\ud{x})=0\,.
$$
Indeed,
$$
\oph_{\ell+1}^\ab\big(x_0,\ldots,x_{\ell-1},
x_\ell+\ophi^\ab(\sigma^{\ell+1}\ud{x})\big)=
\oph_{\ell+1}^\ab\big(x_0,\ldots,x_{\ell-1},
x_\ell^\p+\ophi^\ab(\sigma^{\ell+1}\ud{x}^\p)\big)
$$
Therefore $x^\p_\ell=x_\ell+1$,
$\ophi^\ab(\sigma^{\ell+1}\ud{x})=1$ and
$\ophi^\ab(\sigma^{\ell+1}\ud{x})=0$. Suppose that
$\ophi^\ab(\sigma^k\ud{u})=1$, and apply the above result to
$\sigma^k\ud{u}$ and $\ud{v}$ to get the existence of $m$ with
$\ophi^\ab(\sigma^m\ud{v})=0$. \qed

Let $\ud{u}\in\tA^{\Z_+}$ with $u_0=0$ and
$\ud{u}\preceq\sigma^n\ud{u}$ for all $n\geq 0$. We introduce the
quantity
$$
\widehat{\ud{u}}:=\sup\{\sigma^n\ud{u}\colon n\geq 0\}\,.
$$
We have
$$
\sigma^n\widehat{\ud{u}}\leq\widehat{\ud{u}}\quad\forall n\geq
0\,.
$$
Indeed, if $\widehat{u}$ is periodic, then this is immediate.
Otherwise there exists $n_j$, with $n_j\uparrow\infty$ as
$j\ra\infty$, so that
$\widehat{\ud{u}}=\lim_j\sigma^{n_j}\ud{u}$. By continuity
$$
\sigma^n\widehat{\ud{u}}=
\lim_{j\ra\infty}\sigma^{n+n_j}\ud{u}\leq\widehat{\ud{u}}\,.
$$

\noindent
{\bf Example.\,} We consider the strings $\ud{u}^\p=(01)^\infty$
and $\ud{v}^\p=(110)^\infty$. One can prove that
$\ud{u}^\p=\ud{u}^\ab$ and $\ud{v}^\p=\ud{v}^\ab$ where $\beta$
is the largest solution of
$$
\beta^6-\beta^5-\beta=\beta(\beta^2-\beta+1)(\beta^3-\beta-1)=0
$$
and $\alpha=(1+\beta)^{-1}$. With the notations of  Proposition
\ref{pro3.4ter} we have
$$
\ud{a}=01\quad\ud{a}^\p=00\quad\ud{b}=110\quad\ud{b}^\p=111\,.
$$
Let
$$
\ud{u}:=(00110111)^\infty=(\ud{a}^\p\ud{b}\ud{b}^\prime)^\infty\,.
$$
We have
$$
\widehat{\ud{u}}=(11100110)^\infty=(\ud{b}^\p\ud{a}^\p\ud{b})^\infty\,.
$$
By definition $\ophi^\ab(\sigma\ud{u})=\alpha$. We have
$$
(\ud{b})^\infty\preceq \widehat{\ud{u}}\preceq
\ud{b}^\p(\ud{a})^\infty\,.
$$
From Proposition \ref{pro3.4ter} point 3 and Proposition
\ref{pro3.6} we conclude that
$\log_2\beta=h(\BSigma(\ud{u},\widehat{\ud{u}}))$.

\begin{thm}\label{thm4.1}
Let $k\geq 2$ and let $\ud{u}\in\tA^{\Z_+}$ and
$\ud{v}\in\tA^{\Z_+}$, such that $u_0=0$, $v_0=k-1$ and
\eqref{4.1} holds. If $k=2$ we also assume that
$\sigma\ud{u}\preceq\sigma\ud{v}$. Set
$\log_2\widehat{\beta}:=h(\BSigma(\ud{u},\widehat{\ud{u}}))$. Let
$\bar{\alpha}$ and $\bar{\beta}$ be  defined by the algorithm
of Proposition \ref{pro3.3}. Then\\
1)\, If $\widehat{\beta}<\bar{\beta}$, then $\ud{u}=\ud{u}^\abb$
and $\ud{v}=\ud{v}^\abb $.\\
2)\, If $\widehat{\beta}=\bar{\beta}>1$ and $\ud{u}^{\abb}$ and
$\ud{v}^{\abb}$ are not both periodic, then $\ud{u}=\ud{u}^{\abb}$
and $\ud{v}=\ud{v}^{\abb}$.\\
3)\, If $\widehat{\beta}=\bar{\beta}>1$ and $\ud{u}^{\abb}$ and
$\ud{v}^{\abb}$ are  both periodic, then
$\ud{u}\not=\ud{u}^{\abb}$ and $\ud{v}\not=\ud{v}^{\abb}$.
\end{thm}

\prf Let $\widehat{\beta}<\bar{\beta}$. Suppose  that
$\ud{u}\not=\ud{u}^\abb$ or $\ud{v}\not=\ud{v}^\abb $. By
Proposition \ref{pro3.5} $\ud{u}\not=\ud{u}^\abb$ and
$\ud{v}\not=\ud{v}^\abb $, and there exists $n$ such that
$\ophi^\abb(\sigma^n\ud{u})=1$. Hence
$\ophi^\abb(\widehat{\ud{u}})=1$. If $\bar{\gamma}>0$, then
$\widehat{\ud{u}}_0=v_0=k-1$ whence
$\sigma\widehat{\ud{u}}\preceq\sigma\ud{v}$, so that
$\ophi^\abb(\sigma\widehat{\ud{u}})=\bar{\gamma}$. By Propositions
\ref{pro3.4ter} and \ref{pro3.6} we deduce that
$$
\log_2\widehat{\beta}=h(\BSigma(\ud{u},\widehat{\ud{u}}))
=h(\BSigma(\ud{u},\ud{v}))=\log_2\bar{\beta}\,,
$$
a contradiction. If $\bar{\gamma}=0$, either
$\widehat{\ud{u}}_0=k-1$ and
$\ophi^\abb(\sigma\widehat{\ud{u}})=\bar{\gamma}$, and we get a
contradiction as above, or $\widehat{\ud{u}}_0=k-2$ and
$\ophi^\abb(\sigma\widehat{\ud{u}})=1$. In the latter case, since
$\sigma\widehat{\ud{u}}\preceq \widehat{\ud{u}}$, we conclude that
$\widehat{\ud{u}}_1=k-2$ and
$\ophi^\abb(\sigma^2\widehat{\ud{u}})=1$. Using
$\sigma^n\widehat{\ud{u}}\preceq \widehat{\ud{u}}$ we get
$\widehat{\ud{u}}=(k-2)^\infty=\ud{v}^{\abb}$, so that
$h(\BSigma(\ud{u},\widehat{\ud{u}})) =h(\BSigma(\ud{u},\ud{v}))$,
a contradiction.\\
We prove 2. Suppose for example that $\ud{u}^\abb$ is not
periodic. This implies that $\bar{\alpha}<1$, so that Proposition
\ref{pro3.4} implies that $\ud{u}=\ud{u}^\abb$. We conclude using
Proposition \ref{pro3.5}. Similar proof if $\ud{v}^\abb$ is not
periodic.\\
We prove 3. By Proposition \ref{pro3.5}, $\ud{u}=\ud{u}^{\abb}$ or
$\ud{v}=\ud{v}^{\abb}$ if and only if $\ud{u}=\ud{u}^{\abb}$ and
$\ud{v}=\ud{v}^{\abb}$. Suppose $\ud{u}=\ud{u}^{\abb}$, then
$\ud{u}$ is periodic so that $\widehat{u}=\sigma^p\ud{u}$ for
some $p$. This implies that
$$
\ophi^\abb(\sigma\widehat{\ud{u}})\leq\ophi^\abb(\widehat{\ud{u}})
=\ophi^\abb(\sigma^p\ud{u})<1\,.
$$
by Proposition \ref{pro3.5}. Let $\widehat{\ud{u}}_0\equiv
\widehat{k}-1$. We can apply the algorithm of Proposition
\ref{pro3.3} to the pair $(\ud{u},\widehat{\ud{u}})$ and get two
real numbers $\widetilde{\alpha}$ and $\widetilde{\beta}$ (if
$\widehat{k}=2$, using $\widehat{\beta}>1$ and Theorem
\ref{thm3.1}, we have
$\sigma\ud{u}\preceq\sigma\widehat{\ud{u}}$). Theorem
\ref{thm3.1} implies $\widehat{\beta}=\widetilde{\beta}$, whence
$\widetilde{\beta}=\bar{\beta}$. The map
$\alpha\mapsto\ophi^{\alpha,\bar{\beta}}(\sigma\ud{u})$ is
continuous and decreasing, so that
$\alpha\mapsto\ophi^{\alpha,\bar{\beta}}(\sigma\ud{u})-\alpha$ is
strictly decreasing, whence there exists a unique solution to the
equation $\ophi^{\alpha,\bar{\beta}}(\sigma\ud{u})-\alpha=0$,
which is $\bar{\alpha}=\widetilde{\alpha}$. Therefore
$\ophi^\abb(\sigma\widehat{\ud{u}})<1$ and we must have
$\widehat{k}=k$, whence
$$
\ophi^\abb(\sigma\widehat{\ud{u}})=\bar{\alpha}+\bar{\beta}-k+1=
\ophi^\abb(\sigma\ud{v})\,.
$$
But this implies $\ophi^\abb(\widehat{\ud{u}})=1$, a contradiction
\qed

\begin{thm}\label{thm4.2}
Let $k\geq 2$ and let $\ud{u}\in\tA^{\Z_+}$ and
$\ud{v}\in\tA^{\Z_+}$, such that $u_0=0$, $v_0=k-1$ and
\eqref{4.1} holds. If $k=2$ we also assume that
$\sigma\ud{u}\preceq\sigma\ud{v}$. Let $\bar{\alpha}$ and
$\bar{\beta}$ be  defined by the algorithm of Proposition
\ref{pro3.3}. If $h(\BSigma(\ud{u},\widehat{\ud{u}}))>1$, then
there exists $\ud{u}_*\succeq\widehat{\ud{u}}$ such that
\begin{align*}
\ud{u}_*\prec\ud{v}&\implies \text{$\ud{u}=\ud{u}^\abb$ and
$\ud{v}=\ud{v}^\abb$}\\
\ud{u}_*\succ\ud{v}&\implies \text{$\ud{u}\not=\ud{u}^\abb$ and
$\ud{v}\not=\ud{v}^\abb$}
\end{align*}
\end{thm}

\prf As in the proof of Theorem \ref{thm4.1} we define
$\widetilde{k}$ and, by the algorithm of Proposition \ref{pro3.3}
applied to the pair $(\ud{u},\widehat{\ud{u}})$, two real numbers
$\widetilde{\alpha}$ and $\widetilde{\beta}$. By Theorem
\ref{thm3.1},
$\log_2\widetilde{\beta}=h(\BSigma(\ud{u},\widehat{\ud{u}}))$. We
set
$$
\ud{u}_*:=
\begin{cases}\ud{v}^{\widetilde{\alpha},\widetilde{\beta}}_*&
\text{if $\ud{v}^{\widetilde{\alpha},\widetilde{\beta}}$ is periodic}\\
\ud{v}^{\widetilde{\alpha},\widetilde{\beta}}&\text{if
$\ud{v}^{\widetilde{\alpha},\widetilde{\beta}}$ is not
periodic}\,.
\end{cases}
$$
It is sufficient to show that $\ud{u}_*\prec\ud{v}$ implies
$\bar{\beta}>\widetilde{\beta}$ (see Theorem \ref{thm4.1} point
1). Suppose the contrary, $\bar{\beta}=\widetilde{\beta}$. Then
$$
1=\ophi^{\widetilde{\alpha},\bar{\beta}}(\widehat{\ud{u}})\leq
\ophi^{\widetilde{\alpha},\bar{\beta}}(\ud{v})\,.
$$
We have $\ophi^{\bar{\alpha},\bar{\beta}}(\ud{v})=1$ and for
$\alpha>\bar{\alpha}$, $\ophi^{\alpha,\bar{\beta}}(\ud{v})<1$
(see Lemma \ref{lem3.3.1}). Therefore
$\widetilde{\alpha}\leq\bar{\alpha}$. On the other hand, applying
Corollary \ref{cor3.4} we get $\widetilde{\alpha}\geq\bar{\alpha}$
so that $\widetilde{\alpha}=\bar{\alpha}$ and $\widetilde{k}=k$.
From Propositions \ref{pro3.4bis} or \ref{pro3.4ter} we get
$\ud{v}\preceq \ud{u}_*$, a contradiction.

Suppose that $\ud{u}_*\succ\ud{v}$. We have
$\widehat{u}\preceq\ud{v}\prec \ud{u}_*$, whence
$h(\BSigma(\ud{u},\widehat{\ud{u}}))=h(\BSigma(\ud{u},\ud{u}_*))$
and therefore $\bar{\beta}=\widetilde{\beta}$. As above we show
that $\bar{\alpha}=\widetilde{\alpha}$. Notice that if
$\ud{u}^{\widetilde{\alpha},\widetilde{\beta}}$ is not periodic,
then by Proposition \ref{pro3.4}
$\ud{u}^{\widetilde{\alpha},\widetilde{\beta}}=\ud{u}$. If
$\ud{v}^{\widetilde{\alpha},\widetilde{\beta}}$ is not periodic,
then by Proposition \ref{pro3.4bis}
$\ud{v}^{\widetilde{\alpha},\widetilde{\beta}}=\ud{v}$. If
$\ud{v}^{\widetilde{\alpha},\widetilde{\beta}}$ is  periodic,
then inequalities \eqref{4.1} imply that we must have
$\ud{v}^{\widetilde{\alpha},\widetilde{\beta}}_*\prec \ud{v}$.
Therefore we may have $\ud{u}_*\succ\ud{v}$ and inequalities
\eqref{4.1} only if
$\ud{u}^{\widetilde{\alpha},\widetilde{\beta}}$ and
$\ud{v}^{\widetilde{\alpha},\widetilde{\beta}}$ are periodic.
Suppose that it is the case. If $\ud{u}$ is not periodic, then
using Proposition \ref{pro3.5} the second statement is true. If
$\ud{u}$ is periodic, then $\widehat{\ud{u}}=\sigma^p\ud{u}$ for
some $p$, whence
$\ophi^{\widetilde{\alpha},\widetilde{\beta}}(\sigma^p\ud{u})=1$;
by Proposition \ref{pro3.5}
$\ud{u}\not=\ud{u}^{\widetilde{\alpha},\widetilde{\beta}}$.
 \qed

\end{document}